\documentclass[11pt]{amsart}
\usepackage{amsmath,amssymb}
\usepackage{amsfonts}
\usepackage{amscd}
\usepackage{lscape}
\usepackage{latexsym}
\usepackage{amsfonts,amssymb,stmaryrd,amscd,amsmath,latexsym,amsbsy}

\newcommand{\nc}{\mbox{${\mathbb C}$}}
\newcommand{\nz}{\mbox{${\mathbb Z}$}}

\newcommand{\snz}{\mbox{\scriptsize{${\mathbb Z}$}}}

\newcommand{\gtg}{\mathfrak g}
\newcommand{\gtgl}{\mathfrak gl}
\newcommand{\gtsl}{\mathfrak sl}
\newcommand{\gth}{\mathfrak h}
\newcommand{\gtt}{\mathfrak t}

\newcommand{\tf}{\widetilde{f}}
\newcommand{\te}{\widetilde{e}}
\newcommand{\eps}{\varepsilon}
\newcommand{\vphi}{\varphi}
\newcommand{\cB}{{\mathcal B}}
\newcommand{\cM}{{\mathcal M}}
\newcommand{\cBZ}{{\mathcal BZ}}

\newcommand{\homc}{\mbox{Hom}_{\mbox{\scriptsize{${\mathbb C}$}}}}

\newcommand{\htf}{\widehat{{f}}}
\newcommand{\hte}{\widehat{{e}}}
\newcommand{\heps}{\widehat{\varepsilon}}
\newcommand{\hvphi}{\widehat{\varphi}}

 \newtheorem{thm}{Theorem}[subsection]
 \newtheorem{prop}[thm]{Proposition}
 \newtheorem{lemma}[thm]{Lemma}
 \newtheorem{cor}[thm]{Corollary}
 \newtheorem{defn}[thm]{Definition}

 \newtheorem{rem}{Remark}

\newtheorem{claim}{Claim}[section]

\newenvironment{Ac}%
 {\hspace*{-1.2em}\textbf{Acknowledgment.}\hspace{0.7em}}{}
 {\hspace*{-1.2em}\textbf{Notation.}\hspace{0.7em}}{}

\title[Toward Berenstein-Zelevinsky data in affine type $A$, part II]
{Toward Berenstein-Zelevinsky data in affine type $A$, \\
part II: Explicit description}
\author{Satoshi Naito, Daisuke Sagaki, and Yoshihisa Saito}
\address{Satoshi Naito: Institute of Mathematics, University of Tsukuba, 
Ibaraki 305-8571, Japan.}
\email{naito@math.tsukuba.ac.jp}
\address{Daisuke Sagaki: Institute of Mathematics, University of Tsukuba, 
Ibaraki 305-8571, Japan.}
\email{sagaki@math.tsukuba.ac.jp}
\address{Yoshihisa Saito: Graduate School of Mathematical Sciences, 
University of Tokyo, 3-8-1 Komaba,
Meguro-ku, Tokyo 153-8914, Japan.}
\email{yosihisa@ms.u-tokyo.ac.jp}


\begin{document}
\bigskip
\begin{abstract}
In the present paper, we give an explicit description of the affine analogs
of Berenstein-Zelevinsky data 
constructed in \cite{NSS}, in terms
of certain collections of 
nonnegative integers,  which we 
call Lusztig data of type $A_{l-1}^{(1)}$.
\end{abstract}
\maketitle
\section{Introduction}
\subsection{}
This paper is a continuation of our previous one \cite{NSS}, in which we
introduced Berenstein-Zelevinsky data of type $A_{l-1}^{(1)}$. 

Let us recall the construction of Berenstein-Zelevinsky data of type 
$A_{l-1}^{(1)}$.
We first consider a finite interval $I$ in $\nz$, and 
the finite-dimensional simple Lie algebra $\gtg_I$ of type $A_{|I|}$, where 
$|I|$ denotes the cardinality of $I$; here the set $I$ is thought of 
as the index set of the simple roots of $\gtg_I$. Denote by $\gth_I$ the Cartan 
subalgebra of $\gtg_I$, $W_I$ its Weyl group, and by $\varpi_i^I$, $i\in I$, 
its fundamental weights. The set $\Gamma_I$ of chamber weights is defined to be
$\Gamma_I:=\bigcup_{i\in I}W_I\varpi_i^I\subset\gth_I^*$.
Let ${\bf M}=(M_{\gamma})_{\gamma\in \Gamma_{I}}$ be a Berenstein-Zelevinsky 
datum for $\gtg_I$ in the sense of Kamnitzer (\cite{Kam1}, \cite{Kam2}). 
It is a collection of integers indexed by 
the set $\Gamma_I$, with some additional
conditions called ``{\it edge inequalities}'' and ``{\it tropical Pl\"ucker 
relations}'' (see Definition \ref{defn:fin-BZ}). 
In the general setting of \cite{Kam1} and \cite{Kam2}, the set
of chamber weights is defined to be a subset of $(\gth^{\vee}_I)^*$, where
$\gth^{\vee}_I$ is the Cartan subalgebra
of the (Langlands) dual Lie algebra $\gtg_I^{\vee}$ of $\gtg_I$. 
However, since we focus on type $A$, $(\gth^{\vee}_I)^*$ can be naturally 
identified with $\gth_I^*$. 
We denote by $\cBZ_I$
the set of those Berenstein-Zelevinsky data which satisfy the following 
normalization conditions: $M_{w_0^I\varpi_i^I}=0$
for all $i\in I$, where $w_0^I$ denotes the longest element of $W_I$. 
Kamnitzer showed
that $\cBZ_I$ has a crystal structure, under which the $\cBZ_I$ is isomorphic to
$B(\infty)$ of type $A_{|I|}$ (\cite{Kam1}, \cite{Kam2}).

We note that the family 
$\{\cBZ_I~|~\mbox{$I$ is a finite interval in $\nz$}\}$ forms a 
projective system. The set
$\cBZ_{\snz}$ of Berenstein-Zelevinsky data of type $A_{\infty}$ is defined
to be a kind of projective limit of the projective system above. 
Fix an integer $l\in\nz_{\geq 3}$, and consider a Dynkin diagram 
automorphism $\sigma:\nz\to\nz$ in type $A_{\infty}$ given by $\sigma(p)=p+l$.
There is a naturally induced action  $\sigma$ on $\cBZ_{\snz}$. 
Let $\cBZ_{\snz}^{\sigma}$ denote the fixed point subset of
$\cBZ_{\snz}$ under this action of $\sigma$. Then,
$\cBZ_{\snz}^{\sigma}$ has a crystal structure of type 
$A_{l-1}^{(1)}$, which is induced by the one of $\cBZ_I$.
Moreover, there is a particular connected component 
$\cBZ_{\snz}^{\sigma}({\bf O})$ of 
the crystal $\cBZ_{\snz}^{\sigma}$ (see Subsection 4.3),
which is isomorphic to $B(\infty)$ of type $A_{l-1}^{(1)}$. Thus, we have
a new realization of $B(\infty)$ of type $A_{l-1}^{(1)}$ by  
affine analogs of  Berenstein-Zelevinsky data. 
\subsection{}
In this paper, we consider the set of those Berenstein-Zelevinsky 
data for $\gtg_I$ satisfying another normalization conditions: 
$M_{\varpi_i^I}=0$ for all $i\in I$; it is denoted by $\cBZ_I^{e}$.
In \cite{S}, one of the authors defined a crystal structure on $\cBZ_I^{e}$,
under which the $\cBZ_I^e$ is isomorphic to $B(\infty)$ of type $A_{|I|}$, 
and constructed
an explicit isomorphism $\ast$ between $\cBZ_I$ and $\cBZ_I^{e}$. 

As in the case of $\cBZ_I$, the family 
$\{\cBZ^{e}_I~|~\mbox{$I$ is a finite interval in $\nz$}\}$ forms a 
projective system and we can define $\cBZ_{\snz}^e$ as a kind of projective
limit of this system.
Furthermore, since the map $\ast$ is compatible with the process of taking
this projective limit, we have an isomorphism $\ast$ between 
$\cBZ_{\snz}$ and $\cBZ_{\snz}^{e}$. 
Let us denote by $\bigl(\cBZ_{\snz}^e\bigr)^{\sigma}$ ({\it resp}.,
$\bigl(\cBZ_{\snz}^e\bigr)^{\sigma}({\bf O}^*)$) the image of 
$\cBZ_{\snz}^{\sigma}$ ({\it resp}.,
$\cBZ_{\snz}^{\sigma}({\bf O})$) under the map $\ast$.
Thus, we obtain a crystal $\bigl(\cBZ_{\snz}^e\bigr)^{\sigma}$ of type 
$A_{l-1}^{(1)}$, and its particular connected component  
$\bigl(\cBZ_{\snz}^e\bigr)^{\sigma}({\bf O}^*)$, which is isomorphic to
$B(\infty)$ of type $A_{l-1}^{(1)}$.
\subsection{}
As mentioned in \cite{NSS}, we have not yet found an explicit 
characterization of 
$\cBZ_{\snz}^{\sigma}({\bf O})$ in terms of 
``{\it edge inequalities}'' and ``{\it tropical Pl\"ucker relations}'' 
in type $A_{l-1}^{(1)}$.
Instead, in this paper, we give an explicit description of 
$\cBZ_{\snz}^{\sigma}({\bf O})$ by another approach.

Let us recall the corresponding results in finite type $A$.  
Let $I=[n+1,n+m]$ be a finite interval in $\nz$.
In \cite{S}, one of the authors gave an explicit description of 
$\cBZ_I^{e}$ in terms of Lusztig data associated to $I$. 
Here a Lusztig datum associated to $I$ is a collection
of nonnegative integers indexed by the set 
$\Delta_I^+=\{(i,j)~|~n+1\leq i<j\leq 
n+m+1\}$. Denote by $\cB_I$ the set of all Lusztig data associated to $I$.
It is known that $\cB_I$ parametrizes the canonical basis
(or a PBW basis) of $U_q^-(\gtg_I)$ (\cite{L1}), and has the induced 
crystal structure under which the $\cB_I$ is isomorphic to
$B(\infty)$ of type $A_{m}$ (see \cite{R}, \cite{Sav}, and also \cite{S}).
Motivated by the work \cite{BFZ} of Berenstein-Fomin-Zelevinsky,
one of the authors constructed an explicit isomorphism $\Phi_I$ of 
crystals form $\cB_I$ to $\cBZ_I^{e}$. 

The notion of Lusztig data can be generalized to the affine case. First, we
replace $\Delta_I^+$ by the set $\Delta_{\snz}^+=\{(i,j)~|~i,j\in\nz
\mbox{ with }i<j\}$, and consider a collection 
${\bf a}=(a_{i,j})_{(i,j)\in\Delta_{\snz}^+}$ of nonnegative integers indexed
by $\Delta_{\snz}^+$ such that $a_{i,j}=0$ for every $(i,j)\in\Delta_{\snz}^+$ 
with $j-i\gg 0$, which is called a Lusztig datum associated to $\nz$; 
we denote by $\cB_{\snz}$ the set of all such Lusztig data.
Second, we impose the following two conditions on $\cB_{\snz}$: for each
$(i,j)\in\Delta_{\snz}^+$, $a_{i,j}=a_{i+l,j+l}$ (``cyclicity condition''), 
and there exits at least one $0$ in 
$\{a_{i+s,j+s}~|~0\leq s\leq l-1\}$ (``aperiodicity condition'').
Denote by $\cB_{l-1}^{(1),ap}$ the set of those Lusztig data which satisfy the
conditions above. Then, this set  is naturally identified with 
the set of aperiodic multisegments over $\nz/l\nz$ in the sense of
Lusztig \cite{L2}, which parametrizes the canonical basis of the negative
part of the quantized universal enveloping algebra of type $A_{l-1}^{(1)}$.
The combinatorial description of the induced crystal structure on 
the latter set is given by Leclerc-Thibon-Vasserot \cite{LTV}. 
Hence we can endow $\cB_{l-1}^{(1),ap}$ with a crystal structure
induced by the above identification. 

The main result of this paper is as follows (see Theorem \ref{thm:main} 
for details). We construct a map $\Phi_{\snz}:\cB_{l-1}^{(1),ap}\to 
\bigl(\cBZ_{\snz}^e\bigr)^{\sigma}$ as an affine analog of the map
$\Phi_I:\cB_I\to\cBZ_{I}^e$, and prove that the image of $\cB_{l-1}^{(1),ap}$
under the map $\Phi_{\snz}$ coincides with 
$\bigl(\cBZ_{\snz}^e\bigr)^{\sigma}({\bf O}^*)$. 
Thus, we obtain an isomorphism $\Phi_{\snz}:\cB_{l-1}^{(1),ap}
\overset{\sim}{\to}\bigl(\cBZ_{\snz}^e\bigr)^{\sigma}({\bf O}^*)$. 
Also, by composing the map $\ast$, we obtain an isomorphism 
$\ast\circ\Phi_{\snz}:\cB_{l-1}^{(1),ap}\overset{\sim}{\to}
\cBZ_{\snz}^{\sigma}({\bf O})$. This gives us an explicit 
description of $\cBZ_{\snz}^{\sigma}({\bf O})$, which we
announced in \cite{NSS}.

\medskip
\begin{Ac}
Research of SN is supported by Grant-in-Aid for Scientific 
Research (C) No. 20540006. 
Research of DS is supported by Grant-in-Aid for Young 
Scientists (B) No. 19740004. 
Research of YS is supported by Grant-in-Aid for Scientific 
Research (C) No. 20540009.

The authors are grateful to Professor Saburo Kakei, Professor Yoshiyuki Kimura,
and Professor Yoshihiro Takeyama for valuable discussions. 
\end{Ac}
\section{Review on Berenstein-Zelevinsky data associated to finite 
intervals}
\subsection{Root datum of type $A_m$}
In this subsection, we fix basic notation in type $A_m$,  
following \cite{NSS}; however, for our purpose, we need some changes.
More specifically, in this paper, we realize the root datum of type $A_m$ as
that of $\gtgl_{m+1}(\nc)$, while in \cite{NSS}, we realize it as that of 
$\gtsl_{m+1}(\nc)$.\\

Let $I=[n+1,n+m]$ be a fixed finite interval in $\nz$.
Set $\widetilde{I}:=I\cup\{n+m+1\}$ and consider a finite-dimensional complex 
vector space $\gtt_I$ with basis $\{\epsilon_i\}_{i\in\widetilde{I}}$.
Let $\gtt_I^*:=\homc(\gtt_I,\nc)$ be the dual space of $\gtt_I$. Define
$\Lambda_i^I\in\gtt_I^*$, $i\in I$, by
$$
\langle\epsilon_j,\Lambda_i^I\rangle_I=\left\{\begin{array}{rl}
1 & \mbox{if }j\leq i,\\
0 & \mbox{if }j>i,
\end{array}\right.\quad\mbox{for }i\in I,~j\in \widetilde{I},
$$
where $\langle\cdot,\cdot\rangle_I:\gtt_I\times\gtt^*_I\to\nc$ is the canonical
pairing. If we introduce an element $\kappa_I$ of $\gtt_I^*$ by
$$\langle\epsilon_i,\kappa_I\rangle_I=
1\quad\mbox{for }i\in\widetilde{I},$$
then $\gtt^*_I$ has a basis $\{\Lambda_i^I\}_{i\in I}\cup\{\kappa_I\}$.
For $i\in I$, we define
$$h_i:=\epsilon_i-\epsilon_{i+1}\in\gtt_I\quad\mbox{and}\quad
\alpha_i^I:=-\Lambda^I_{i-1}+2\Lambda_i^I-\Lambda_{i+1}^I\in\gtt^*_I.$$
Here we set $\Lambda^I_{n-1}=0$ and $\Lambda^I_{n+m+1}=\kappa_I$ 
by convention. Then, 
$\{\alpha_i^I\}_{i\in I}\cup\{\kappa_I\}$ forms another basis of $\gtt^*_I$:
$$\gtt^*_I=\left(\mathop{\oplus}_{i\in I}\nc\alpha_i^I\right)\oplus
\nc\kappa_I.\eqno{(2.1.1)}$$
Also, the matrix 
$A_I=(a_{ij})_{i,j\in I}:=
\left(\langle h_i,\alpha_j^I\rangle_I\right)_{i,j\in I}$ is the 
Cartan matrix of type $A_m$ indexed by $I$. Namely, for $i,j\in I$, 
$$a_{ij}=\langle h_j,\alpha_i^I\rangle_I=\left\{\begin{array}{rl}
2 & \mbox{if }j=i,\\
-1 & \mbox{if }j=i\pm 1,\\
0 & \mbox{otherwise}.
\end{array}\right.$$
Let $\gth_I$ be the subspace of $\gtt_I$ (of codimension $1$) spanned by 
the set $\Pi^{\vee}_I:=\{h_i\}_{i\in I}$. Then there is a canonical projection 
$\pi_I:\gtt_I^*\to \gth_I^*:=\homc(\gth_I,\nc)$, for which 
$\mbox{Ker}(\pi_I)=\nc\kappa_I$ holds.
By the splitting (2.1.1), we can naturally identify $\gth_I^*$ with
the subspace $(\gtt^*_I)^0$ of $\gtt^*_I$ spanned by 
$\{\alpha_i^I\}_{i\in I}$. Let us denote the induced embedding 
$\iota_I:\gth_I^*\hookrightarrow \gtt_I^*$.
Set $\Pi_I:=\{\pi_I(\alpha_i^I)\}_{i\in I}~(\subset \gth_I^*)$. 
Thus, $(A_I; \Pi_I,\Pi^{\vee}_I,\gth_I^*,\gth_I)$ is a root datum 
of type $A_m$. Under the above identification, we can write
$\Pi_I=\{\alpha_i^I\}_{i\in I}~\subset\gtt_I^*$.

Let $\gtg_I$ be the complex finite-dimensional simple Lie algebra with Cartan
matrix $A_I$ and Cartan subalgebra $\gth_I$. Then, $\Pi_I=\{h_i\}_{i\in I}$ and 
$\Pi_I^{\vee}=\{\pi_I(\alpha_i)\}_{i\in I}$ are the set of simple
coroots and roots of $\gtg_I$, respectively.  
Also, $\varpi_i^I:=\pi_I(\Lambda_i^I)\in\gth_I^*$, $i\in I$, are the fundamental
weights for $\gtg_I$. \\

We define a linear automorphism $\sigma_i$, $i\in I$,  of $\gtt_I^*$ by
$$\sigma_i(\lambda):=\lambda-\langle h_i,\lambda\rangle_I\alpha_i^I
\quad \mbox{for }\lambda\in\gtt_I^*,$$
and the linear automorphism $s_i$, $i\in I$ of $\gth_I^*$ by
$$s_i:=\pi_I\circ \sigma_i\circ \iota_I.$$
Then we have
$$s_i(\mu)=\mu-\langle h_i,\mu\rangle \pi_I(\alpha_i^I)\quad\mbox{for }
\mu\in \gth_I^*,$$
where $\langle\cdot,\cdot\rangle:\gth_I\times\gth_I^*\to \nc$ is the canonical
pairing. 
The group $W_{\gth_I^*}:=\langle s_i~|~i\in I\rangle~\bigl(\subset\mbox{Aut}
(\gth_I^*)\bigr)$ is thought of as the Weyl group of $\gtg_I$. 
Because $W_{\gth_I^*}$ is isomorphic to $W_{\gtt_I^*}:=
\langle \sigma_i~|~i\in I\rangle~\bigl(\subset\mbox{Aut}(\gtt_I^*)\bigr)$ via
the map $s_i\mapsto \sigma_i$, we write $W_I=W_{\gtt_I^*}$ for the remainder
of this paper.
\subsection{Chamber weights}
Consider the set $\{w\varpi_i^I~|~w\in W_{\gth_I^*},~i\in I\}~
\bigl(\subset\gth_I^*\bigr)$. In \cite{NSS},
this set is denoted by $\Gamma_I$, and an element of it is called a chamber 
weight. However, in this paper, we change the notation. Namely, we set
$$\Gamma_I:=\{w\Lambda_i^I~|~w\in W_I,~i\in I\}~
\bigl(\subset\gtt_I^*\bigr),$$
and call an element of this set a chamber weight. Since there is a 
bijection $\Gamma_I\overset{\sim}{\to}\{w\varpi_i^I~|~w\in W_{\gth_I^*},
~i\in I\}$ induced by $\pi_I:\gtt_I^*\to\gth_I^*$, these two sets are
naturally identified. Thus, the set of chamber weights in the sense
of \cite{NSS} is equal to $\pi_I(\Gamma_I)$.

Let $w_0^I$ denote the longest element of $W_{\gth_I^*}$, and $\omega_I:I\to I$ 
an involution defined by $n+i\mapsto n+m-i+1$, $1\leq i\leq m$. 
As in \cite{NSS},
the following two equalities are verified in $\gth_I^*$ and in $W_{\gth_I^*}$, 
respectively: for $i\in I$, 
$$w_0^I(\varpi_i^I)=-\varpi_{\omega_I(i)}^I,\quad w_0^Is_{\omega_I(i)}=s_iw_0^I.
\eqno{(2.2.1)}$$
Let $\gamma'\in \pi_I(\Gamma_I)$. Recall that it can be written as 
$\gamma'=w\varpi^I_i$ for some $w\in W_{\gth_I^*}$ and $i\in I$. 
Because of (2.2.1), we have
$$-\gamma'=ww_0^I\varpi_{\omega_I(i)}\eqno{(2.2.2)}$$
and
$$\{-w\varpi_i^I~|~w\in W_{\gth_I^*},~i\in I\}=
\{w\varpi_i^I~|~w\in W_{\gth_I^*},~i\in I\}.$$

Now we consider analogs of (2.2.1) for 
$\gtt_I^*$. As in the case above, we have
$$w_0^I\sigma_i=\sigma_{\omega_I(i)}w_0^I,\eqno{(2.2.3)}$$
where we also denote by $w_0^I$ the longest element of $W_I=W_{\gtt_I^*}.$
However, an obvious analog of the first equality of (2.2.1) fails. 
In fact, we have
$$w_0^I(\Lambda_i^I)=\kappa_I-\Lambda_{\omega_I(i)}^I.$$
As a consequence, $\Gamma_I$ does not coincide with $-\Gamma_I$. Still, 
we can define an analog of the map $\gamma'\mapsto -\gamma'$ as follows. For 
$\gamma=w\Lambda_i^I$ with $w\in W_I,i\in I$, we set
$$\gamma^c:=ww_0^I\Lambda_{\omega_I(i)}.$$
Then, it is easy to see that the map $c:\gamma\mapsto\gamma^c$
is well-defined, and gives an involution of $\Gamma_I$. In this
paper, this map plays an important role.\\

A set of integers ${\bf k}=\{k_{1}<k_{2}<\cdots<k_u\}\subset 
\widetilde{I}$ is called a Maya diagram associated to $I$. We denote by 
$\cM_I$ the set of all 
Maya diagrams associated to $I$, and set 
$\cM_I^{\times}:=\cM_I\setminus \{\phi,\widetilde{I}\}$. For
a given chamber weight $\gamma\in \Gamma_I$, we define the associated Maya
diagram ${\bf k}(\gamma)$ by
$${\bf k}(\gamma):=\{k\in \widetilde{I}~|~\langle \epsilon_k,\gamma
\rangle_I=1\}.$$
Note that $\langle \epsilon_k,\gamma\rangle_I=0\mbox{ or }1~
\mbox{for each }k\in\widetilde{I}\mbox{ and }\gamma\in\Gamma_I$
by the definitions.
Since the map $\gamma\mapsto {\bf k}(\gamma)$ gives a bijection from
$\Gamma_I$ to $\cM_I^{\times}$, we can identify these two sets.  
 
Under the identification above, the involution $c$ is described as follows
(see \cite{S}). For a Maya diagram ${\bf k}\in \cM_I^{\times}$, define
${\bf k}^c:=\widetilde{I}\setminus {\bf k}$. Then the next lemma follows easily.
\begin{lemma}
$${\bf k}(\gamma^c)={\bf k}(\gamma)^c\quad(\gamma\in\Gamma_I).$$
\end{lemma}
\subsection{BZ data associated to $I$}
Let
${\bf M}=(M_{\bf k})_{{\bf k}\in\cM_I^{\times}}$ be a collection of integers
indexed by $\cM_I^{\times}$. For each ${\bf k}\in\cM_I^{\times}$, we call
$M_{\bf k}$ the ${\bf k}$-component of ${\bf M}$ and denote
it by $({\bf M})_{\bf k}$. Under the identification 
$\Gamma_I\cong \cM_I^{\times}$, Berenstein-Zelevinsky data
are defined as follows (see \cite{S}).
\begin{defn}\label{defn:fin-BZ}
A collection ${\bf M}=(M_{\bf k})_{{\bf k}\in\cM_I^{\times}}$ of integers 
indexed by $\cM_I^{\times}$ is
called a Berenstein-Zelevinsky {\rm (}BZ for short{\rm )} datum associated to 
an interval $I$ if it satisfies the following conditions:
\vskip 1mm
\noindent
{\rm (BZ-1)} for all indices $i\ne j$ in 
$\widetilde{I}$ and all 
${\bf k}\in \cM_n$ 
with ${\bf k}\cap\{i,j\}=\phi$, 
$$M_{{\bf k}\cup\{i\}}+M_{{\bf k}\cup\{j\}}\leq M_{{\bf k}\cup\{i,j\}}
+M_{\bf k};$$
\noindent
{\rm (BZ-2)} for all indices 
$i<j<k$ in $\widetilde{I}$ and all 
${\bf k}\in \cM_n$ with ${\bf k}\cap\{i,j,k\}=\phi$, 
$$M_{{\bf k}\cup\{i,k\}}+M_{{\bf k}\cup\{j\}}=\mbox{\rm min}\left\{
M_{{\bf k}\cup\{i,j\}}+M_{{\bf k}\cup\{k\}},~
M_{{\bf k}\cup\{j,k\}}+M_{{\bf k}\cup\{i\}}\right\}.$$
Here we set $M_{\phi}=M_{\widetilde{I}}=0$ by convention.  
\end{defn}

\begin{defn}
A BZ datum ${\bf M}=(M_{\bf k})_{{\bf k}\in\cM_I^{\times}}$ is called
a $w_0$-BZ {\rm (}resp., $e$-BZ {\rm )} datum if it satisfies the 
following normalization conditions:
\vskip 1mm
\noindent
{\rm (BZ-0)} for every $i\in I$,
$$M_{{\bf k}(\Lambda_i^I)^c}=0\qquad(\mbox{resp., }M_{{\bf k}(\Lambda_i^I)}=0).$$
\end{defn}

We denote by $\mathcal{BZ}_I$ ({\it resp}., $\mathcal{BZ}_I^{e}$)
the set of all $w_0$-BZ ({\it resp}., $e$-BZ) data. Consider a 
map $\ast:\mathcal{BZ}_I\to \mathcal{BZ}_I^{e}$ by
$$({\bf M}^*)_{\bf k}:=({\bf M})_{{\bf k}^c}\quad
({\bf M}\in \mathcal{BZ}_I).$$
By \cite{S}, this is a well-defined bijection; the inverse of it is also 
denoted by $\ast$.\\

Let us discuss crystal structures on BZ data. First,  we define a crystal
structure on $\cBZ_I$, following Kamnitzer \cite{Kam2} 
(see also \cite{S}).
For ${\bf M}=(M_{\bf k})_{{\bf k}\in\cM_I^{\times}}\in \cBZ_I$,
we define the weight $\mbox{wt}({\bf M})$ of ${\bf M}$ by
$$\mbox{wt}({\bf M}):=\sum_{i\in I}M_{{\bf k}(\Lambda_i^I)}\alpha_i^I.$$
For $i\in I$, we set
$$\eps_i({\bf M}):=-\left(M_{{\bf k}(\Lambda_i^I)}
+M_{{\bf k}(\sigma_i\Lambda_i^I)}-M_{{\bf k}(\Lambda_{i-1}^I)}
-M_{{\bf k}(\Lambda_{i+1}^I)}\right),$$
$$
\vphi_i({\bf M}):=\eps_i({\bf M})
+\langle h_i,\mbox{wt}({\bf M})\rangle_I.
$$
In order to define the action of Kashiwara operators $\te_i$ and $\tf_i$, 
$i\in I$, we recall the following result, due to Kamnitzer:
\begin{prop}[\cite{Kam2}]
Let ${\bf M}=(M_{\bf k})\in \mathcal{BZ}_I$ be a $w_0$-BZ 
datum associated to $I$.
\vskip 1mm
\noindent
{\rm (1)} If $\eps_i({\bf M})>0$, then there exists a unique $w_0$-BZ 
datum, denoted by 
$\te_i{\bf M}$, such that
\begin{itemize}
\item[(i)] $(\te_i{\bf M})_{{\bf k}(\Lambda_i^I)}=M_{{\bf k}(\Lambda_i^I)}+1$,
\item[(ii)] $(\te_i{\bf M})_{\bf k}=M_{\bf k}$ for all 
${\bf k}\in \cM_I^{\times}\setminus \cM_I^{\times}(i)$.
\end{itemize}
Here, $\cM_I^{\times}(i)=\left\{{\bf k}\in \cM_I^{\times}~|~i\in {\bf k}
\mbox{ and }i+1\not\in {\bf k}\right\}\subset \cM_I^{\times}$.
\vskip 1mm
\noindent
{\rm (2)} There exists a unique $w_0$-BZ datum, denoted by 
$\tf_i{\bf M}$, such that
\begin{itemize}
\item[(iii)] $(\tf_i{\bf M})_{{\bf k}(\Lambda_i^I)}=M_{{\bf k}(\Lambda_i^I)}-1$,
\item[(iv)] $(\tf_i{\bf M})_{\bf k}=M_{\bf k}$ for all 
${\bf k}\in \cM_n^{\times}\setminus \cM_n^{\times}(i)$.
\end{itemize}
\end{prop}
If $\eps_i({\bf M})=0$, then we set $\te_i{\bf M}={0}$. 
\begin{thm}[\cite{Kam2}]\label{thm:BZ-crys}
The set $\mathcal{BZ}_I$, equipped with the maps 
$\mbox{\em wt},\eps_i,\vphi_i,\te_i,\tf_i$, is a crystal
which is isomorphic to 
$\left(B(\infty);\mbox{\rm wt},\eps_i,\vphi_i,\te_i,\tf_i\right)$. 
\end{thm}
\vskip 3mm
Second, we introduce a crystal structure on $\cBZ_I^e$.
For ${\bf M}=(M_{\bf k})_{{\bf k}\in\cM_I^{\times}}\in \cBZ_I^{e}$,
we set
$$\mbox{wt}({\bf M}):=\mbox{wt}({\bf M}^*),\quad
\eps_i^*({\bf M}):=\eps_i({\bf M}^*),\quad
\vphi_i^*({\bf M}):=\vphi_i({\bf M}^*).$$
\begin{cor}[\cite{S}]\label{cor:ast-crys}
{\rm (1)} Let ${\bf M}=(M_{\bf k})\in \mathcal{BZ}^{e}_I$ be an $e$-BZ 
datum associated to $I$. 
If $\eps_i^*({\bf M})>0$, then there exists a unique $e$-BZ datum,
denoted by $\te_i^*{\bf M}$, such that
\begin{itemize}
\item[(i)] $(\te_i^*{\bf M})_{{\bf k}(\Lambda_i^I)^c}=M_{{\bf k}(\Lambda_i^I)^c}+1$,
\item[(ii)] $(\te_i^*{\bf M})_{\bf k}=M_{\bf k}$ for all 
${\bf k}\in \cM_I^{\times}\setminus \cM_I^{\times}(i)^{\ast}$. 
\end{itemize}
Here, $\cM_I^{\times}(i)^*=\left\{{\bf k}\in \cM_I^{\times}~|~
i\not\in {\bf k}\mbox{ and }i+1\in {\bf k}\right\}\subset \cM_I^{\times}$.
\vskip 1mm
\noindent
{\rm (2)} There exists a unique $e$-BZ datum, denoted by 
$\tf_i^*{\bf M}$,  such that
\begin{itemize}
\item[(iii)] $(\tf_i^*{\bf M})_{{\bf k}(\Lambda_i^I)^c}=M_{{\bf k}(\Lambda_i^I)^c}-1$,
\item[(iv)] $(\tf_i^*{\bf M})_{\bf k}=M_{\bf k}$ for all 
${\bf k}\in \cM_I^{\times}\setminus \cM_I^{\times}(i)^*$.
\end{itemize}
\vskip 1mm
\noindent
{\rm (3)} If $\eps_i^*({\bf M})=0$, then we set $\te_i^*{\bf M}:={0}$.
The following equalities hold:
$$\te_i^*{\bf M}=(\te_i({\bf M}^*))^*,\quad
\tf_i^*{\bf M}=(\tf_i({\bf M}^*))^*.$$
Here it is understood that ${0}^*={0}.$
\vskip 1mm
\noindent
{\rm (4)} 
The set $\mathcal{BZ}^e_I$, equipped with the maps 
$\mbox{\em wt},\eps_i^*,\vphi_i^*,\te_i^*,\tf_i^*$, is a crystal
which is isomorphic to 
$\left(B(\infty);\mbox{\rm wt},\eps_i^*,\vphi_i^*,\te_i^*,\tf_i^*\right)$. 
Moreover, the bijection $\ast:
\cBZ_I\overset{\sim}{\to}\cBZ_I^e$ gives an isomorphism of crystals
form $\left(\cBZ_I;\mbox{\rm wt},\eps_i,\vphi_i,\te_i,\tf_i\right)$
to $\left(\cBZ_I^{e};\mbox{\rm wt},\eps_i^*,\vphi_i^*,\te_i^*,\tf_i^*\right)$.
\end{cor}
\subsection{Lusztig data associated to $I$}
Recall that $I=[n+1,n+m]$. Let $\Delta^+_I=\{(i,j)~|~i,j\in \widetilde{I}
\mbox{ with }i<j\}$, and consider
$$\cB_I:=\left\{{\bf a}=(a_{i,j})_{(i,j)\in\Delta^+_I}~|~
a_{i,j}\in\nz_{\geq 0}\mbox{ for all }(i,j)\in \Delta^+_I\right\},$$
which is the set of all $m(m+1)/2$-tuples of nonnegative integers 
indexed by $\Delta^+_I$. In this paper, an element of $\cB_I$ is 
called a Lusztig datum associated to $I$.

We define two crystal structures on $\cB_I$, following \cite{S}. 
For ${\bf a}\in \cB_I$, 
define the weight $\mbox{wt}({\bf a})$
of ${\bf a}$ by 
$$\mbox{wt}({\bf a}):=-\sum_{i\in I}r_i({\bf a})\alpha_i^I,\quad\mbox{where}\quad
r_i({\bf a}):=\sum_{s=n+1}^i\sum_{t=i+1}^{n+m+1}a_{s,t}\quad (i\in I).$$
For $i\in I$, set
$$A^{(i)}_k({\bf a}):=\sum_{s=n+1}^k(a_{s,i+1}-a_{s-1,i})\quad 
(n+1\leq k\leq i),$$
$$A^{\ast (i)}_k({\bf a}):=\sum_{t=k+1}^{n+m+1}(a_{i,t}-a_{i+1,t+1})\quad 
(i\leq k\leq n+m).$$
Here we set $a_{n,i}=0$ and $a_{i+1,n+m+2}=0$ by convention. Next, define 
$$\eps_i({\bf a}):=
\mbox{max}\left\{A_{n+1}^{(i)}({\bf a}),\ldots,A_i^{(i)}({\bf a})\right\},
\quad\vphi_i({\bf a}) = \eps_i({\bf a})+\langle h_i,\mbox{wt}({\bf a})
\rangle,$$
$$\eps_i^*({\bf a}):=\mbox{max}\left\{A^{\ast (i)}_i({\bf a}),\ldots,
A^{\ast (i)}_{n+m}({\bf a})\right\},\quad
\varphi_i^*({\bf a})=\eps_i^*({\bf a})+\langle h_i,\mbox{wt}({\bf a})
\rangle.$$
We set
$$k_e:=\mbox{min}\left\{n+1\leq k\leq i~\left|~ \eps_i({\bf a})
=A_k^{(i)}({\bf a})\right.\right\},$$
$$k_f:=\mbox{max}\left\{n+1\leq k\leq i~\left|~ \eps_i({\bf a})
=A_k^{(i)}({\bf a})\right.\right\},$$
$$k_e^*:=\mbox{max}\left\{i\leq k\leq n+m~\left|~ \eps_i^*({\bf a})
=A_k^{\ast (i)}({\bf a})\right.\right\},$$
$$k_f^*:=\mbox{min}\left\{i\leq k\leq n+m~\left|~ \eps_i^*({\bf a})
=A_k^{\ast (i)}({\bf a})\right.\right\}.$$
For a given ${\bf a}\in \cB_I$, we introduce the following $m(m+1)/2$-tuples 
of integers 
${\bf a}^{(\sharp)}=\left(a_{s,t}^{(\sharp)}\right)~(\sharp=1,2,3,4)$ by 
\begin{align*}
{a}^{(1)}_{s,t}&:=\left\{\begin{array}{ll}
a_{k_e,i}+1 & \mbox{if }s=k_e,~t=i,\\
a_{k_e,i+1}-1 & \mbox{if }s=k_e,~t=i+1,\\
a_{s,t} & \mbox{otherwise}.
\end{array}\right.\\
{a}_{s,t}^{(2)}&:=\left\{\begin{array}{ll}
a_{k_f,i}-1 & \mbox{if }s=k_f,~t=i,\\
a_{k_f,i+1}+1 & \mbox{if }s=k_f,~t=i+1,\\
a_{s,t} & \mbox{otherwise},
\end{array}\right.\\
{a}_{s,t}^{(3)}&:=\left\{\begin{array}{ll}
a_{i,k_e^*+1}-1 & \mbox{if }s=i,~t=k_e^*+1,\\
a_{i+1,k_e^*+1}+1 & \mbox{if }s=i+1,~t=k_e^*+1,\\
a_{s,t} & \mbox{otherwise}.
\end{array}\right.\\
{a}^{(4)}_{s,t}&:=\left\{\begin{array}{ll}
a_{i,k_f^*+1}+1 & \mbox{if }s=i,~t=k_f^*+1,\\
a_{i+1,k_f^*+1}-1 & \mbox{if }s=i+1,~t=k_f^*+1,\\
a_{s,t} & \mbox{otherwise}.
\end{array}\right.
\end{align*}
Finally, we define Kashiwara operators on $\cB_I$ by:
$$\te_i{\bf a}=\left\{\begin{array}{ll}
{\bf 0}&\mbox{if }\eps_i({\bf a})=0,\\
{\bf a}^{(1)}&\mbox{if }\eps_i({\bf a})>0,
\end{array}\right.\quad \tf_i{\bf a}={\bf a}^{(2)},$$
$$\te_i^*{\bf a}=\left\{\begin{array}{ll}
{\bf 0}&\mbox{if }\eps_i^*({\bf a})=0,\\
{\bf a}^{(3)}&\mbox{if }\eps_i^*({\bf a})>0,
\end{array}\right.\quad \tf_i^*{\bf a}={\bf a}^{(4)}.$$

\begin{prop}[\cite{R},\cite{S}]
Each of $\left(\cB_I;\mbox{\rm wt},\eps_i,\vphi_i,\te_i,\tf_i\right)$
and $\left(\cB_I;\mbox{\rm wt},\eps_i^*,\vphi_i^*,\te_i^*,\tf_i^*\right)$ is 
a crystal which is isomorphic to $B(\infty)$.
\end{prop} 

This proposition tells us that $\cB_I$ has two different crystal structures, 
which are isomorphic to each other.
In this paper, we call the first one the ordinary crystal structure on 
$\cB_I$; the second one is called the $\ast$-crystal structure on $\cB_I$.

\begin{defn}[\cite{BFZ}]
Let ${\bf k}=\{k_{n+1}<k_{n+2}<\cdots<k_{n+u}\}\in\cM_I^{\times}$. 
For such a ${\bf k}$,
we define a ${\bf k}$-tableau as an uppertriangular matrix 
$C=(c_{p,q})_{n+1\leq p\leq q\leq n+u}$ with integer entries satisfying
$$c_{p,p}=k_p\qquad (n+1\leq p\leq n+u),$$
and the {\rm (}usual{\rm )} monotonicity conditions for semi-standard tableaux:
$$c_{p,q}\leq c_{p,q+1},\qquad c_{p,q}<c_{p+1,q}.$$
\end{defn}
For a given Lusztig datum ${\bf a}=(a_{i,j})\in\cB_I$, let ${\bf M}({\bf a})=
(M_{\bf k}({\bf a}))_{{\bf k}\in \cM_I^{\times}}$ be a collection of integers
defined by
$$M_{\bf k}({\bf a}):=-\sum_{j=n+1}^{n+u}\sum_{i=n+1}^{k_j-1}a_{i,k_j}+
\mbox{min}\left\{\left.\sum_{n+1\leq p<q\leq n+u}a_{c_{p,q},c_{p,q}+(q-p)}
~\right|~\begin{array}{c}C=(c_{p,q})\mbox{ is }\\ 
\mbox{a ${\bf k}$-tableau}\end{array}\right\}.\eqno{(2.4.1)}$$
We denote the map ${\bf a}\mapsto {\bf M}({\bf a})$ by $\Phi_I$.
\begin{thm}[\cite{BFZ}, \cite{S}]\label{thm:finA}
For each ${\bf a}\in\cB_I$, $\Phi_I({\bf a})={\bf M}({\bf a})$ is an 
$e$-BZ datum. Moreover, the map $\Phi_I:\cB_I\to \mathcal{BZ}^{e}_I$ is 
a bijection, which induces an isomorphism of crystals
$$\left(\cB_I;\mbox{\rm wt},\eps_i^*,\vphi_i^*,\te_i^*,\tf_i^*\right)
\overset{\Phi_I}{\longrightarrow}
\left(\cBZ_I^{e};\mbox{\rm wt},\eps_i^*,\vphi_i^*,\te_i^*,\tf_i^*\right).$$
\end{thm}
Consider a composition of bijections 
$\cB_I\overset{\Phi_I}{\longrightarrow}\cBZ_I^e
\overset{\ast}{\longrightarrow}\cBZ_I$. By Corollary \ref{cor:ast-crys}
and Theorem \ref{thm:finA}, we obtain the following corollary.
\begin{cor}\label{cor:ast-isom}
The bijection $\ast\circ\Phi_I:\cB_I\to\cBZ_I$
gives an isomorphism of crystals 
$$\left(\cB_I;\mbox{\rm wt},\eps_i^*,\vphi_i^*,\te_i^*,\tf_i^*\right)
\overset{\ast\circ\Phi_I}{\longrightarrow}
\left(\cBZ_I;\mbox{\rm wt},\eps_i,\vphi_i,\te_i,\tf_i\right).$$
\end{cor}
\section{Berenstein-Zelevinsky data associated to $\nz$}
\subsection{Root datum of $\gtgl_{\infty}(\nc)$}
In this subsection, we introduce basic notation in type $A_{\infty}$,
following \cite{NSS}; however, for our purpose, we need some changes 
as in the finite interval case.
More specifically, we need to replace $\gtsl_{\infty}(\nc)$ by
$\gtgl_{\infty}(\nc)$.\\  

Let $\gtt$ be a vector space over $\nc$ with basis 
$\{\epsilon_i\}_{i\in \snz}$. Set
$$h_i:=\epsilon_i-\epsilon_{i+1}\quad\mbox{for }i\in\nz,$$
and consider the subspace $\gth$ of $\gtt$ spanned by 
$\Pi=\{h_i\}_{i\in\snz}$. Note that $\gth$ is a subspace
of $\gtt$ of codimension $1$. 

Let $\gtt^*=\homc(\gtt,\nc)$ ({\it resp}., $\gth^*=\homc(\gth,\nc)$) be the dual
space of $\gtt$ ({\it resp}., $\gth$). Then there is a canonical projection
$\pi:\gtt^*\to\gth^*$, whose kernel is spanned by an element $\kappa$ of 
$\gtt^*$ defined by $\kappa(\epsilon_i)=1$ for all $i\in\nz$. We remark that
there is a splitting of $\gtt^*$:
$$\gtt^*\cong \gth^*\oplus\nc\kappa.\eqno{(3.1.1)}\label{split}$$ 
\vskip 3mm
Define $\Lambda_i,~\Lambda_i^c\in \gtt^*$, $i\in \nz$, by
$$
\langle\epsilon_j,\Lambda_i\rangle_{\snz}:=\left\{\begin{array}{rl}
1 & \mbox{if }j\leq i,\\
0 & \mbox{if }j>i,
\end{array}\right.\qquad
\langle\epsilon_j,\Lambda_i^c\rangle_{\snz}:=\left\{\begin{array}{rl}
0 & \mbox{if }j\leq i,\\
1 & \mbox{if }j>i,
\end{array}\right.
$$
where $\langle\cdot,\cdot\rangle_{\snz}:\gtt\times\gtt^*\to\nc$ 
is the canonical
pairing. Then we have
$$
\langle h_j,\Lambda_i\rangle_{\snz}=\delta_{ij},~
\langle\epsilon_0,\Lambda_i\rangle_{\snz}=\left\{\begin{array}{rl}
1 & \mbox{if }i\geq 0,\\
0 & \mbox{if }i<0,
\end{array}\right.~
\langle h_j,\Lambda_i^c\rangle_{\snz}=-\delta_{ij},~
\langle \epsilon_0,\Lambda_i^c\rangle_{\snz}=\left\{\begin{array}{rl}
0 & \mbox{if }i\geq 0,\\
1 & \mbox{if }i<0.
\end{array}\right.
$$
These formulas enable us to deduce the following lemma:
\begin{lemma}\label{lemma:gltosl}
For each $i\in \nz$, we have $\pi(\Lambda^c_i)=-\pi(\Lambda_i)$.
\end{lemma}
Set
$$\alpha_i:=-\Lambda_{i-1}+2\Lambda_i-\Lambda_{i+1}\quad \mbox{for }i\in\nz.$$
By the definitions, we have
$$\langle h_j,\alpha_i\rangle_{\snz}=\left\{\begin{array}{rl}
2 & \mbox{if }j=i,\\
-1 & \mbox{if }j=i\pm 1,\\
0 & \mbox{otherwise},
\end{array}\right.
\qquad \langle\epsilon_j,\alpha_i\rangle_{\snz}=
\left\{\begin{array}{rl}
1 & \mbox{if }i=j,\\
-1 & \mbox{if }i=j-1,\\
0 & \mbox{otherwise}.
\end{array}\right.
$$
Let $\langle\cdot,\cdot\rangle:\gth\times\gth^*\to\nc$ be the canonical
pairing. Then, $\langle h_j,\pi(\alpha_i)\rangle=\langle h_j,
\alpha_i\rangle_{\snz}$, and $\Pi^{\vee}:=\{\pi(\alpha_i)\}_{i\in\nz}$ is 
a linearly independent subset
of $\gth^*$. In other words, $(A_{\snz},\Pi,\Pi^{\vee},\gth^*,\gth)$ is a root
datum of type $A_{\infty}$ in the sense of \cite{NSS}. Here, 
$A_{\snz}=(a_{ij})_{i,j\in\snz}$ is the Cartan matrix of type $A_{\infty}$, whose
entries $a_{ij}$ are defined by $a_{ij}:=\langle h_j,\pi(\alpha_i)\rangle$
for $i,j\in\nz$.\\

For each $i\in\nz$, Define $\sigma_i\in\mbox{Aut}(\gtt^*)$ by
$$\sigma_i(\lambda):=\lambda-\langle h_i,\lambda\rangle_{\snz} \alpha_i\quad
\mbox{for }\lambda\in \gtt^*.$$
Similarly, define $\sigma_i\in \mbox{Aut}(\gtt)$ by
$$\sigma_i(t):=t-\langle t,\alpha_i\rangle_{\snz} h_i\quad
\mbox{for }t\in \gtt.$$
Let $\iota:\gth^*\hookrightarrow \gtt^*$ be the embedding induced by the
splitting (3.1.1). If we define $s_i\in\mbox{Aut}(\gth^*)$ by
$$s_i:=\pi\circ\sigma_i\circ\iota,$$
then we have
$$s_i(\mu)=\mu-\langle h_i,\mu\rangle\pi(\alpha_i)\quad
\mbox{for }\mu\in\gth^*.$$
Let $W_{\snz}:=\langle s_i~|~i\in\nz\rangle\subset \mbox{Aut}(\gth^*)$ 
be the Weyl group of type $A_{\infty}$. Consider the subgroup of 
$\mbox{Aut}(\gtt^*)$ generated by $\sigma_i,~i\in \nz$. It is easy to see 
that this group is isomorphic to $W_{\snz}$. For this reason, we also denote it
by $W_{\snz}$.
\subsection{Chamber weights associated to $\nz$ and Maya diagrams}
Set
$$\Xi_{\snz}=\{w\Lambda_i~|~w\in W_{\snz},~i\in\nz\}\quad\mbox{and}\quad
\Gamma_{\snz}=\{w\Lambda_i^c~|~w\in W_{\snz},~i\in\nz\}.$$
An element of $\Xi_{\snz}$ ({\it resp}., $\Gamma_{\snz}$) is called a chamber
weight ({\it resp}., dual chamber weight) associated to $\nz$.
\begin{defn}
{\rm (1)} For a given integer $r\in\nz$, a subset ${\bf k}$ of $\nz$ is 
called a Maya diagram of charge $r$ if it satisfies the following condition: 
there exist nonnegative integers $p$ and $q$ such that
$$\nz_{\leq r-p}\subset {\bf k}\subset \nz_{\leq r+q},\quad
|{\bf k}\cap\nz_{>r-p}|=p,\eqno{(3.2.1)}$$
where $|{\bf k}\cap\nz_{>r-p}|$ denotes the cardinality of the finite set 
${\bf k}\cap\nz_{>r-p}$.
We denote by $\cM_{\snz}^{(r)}$ the set of all Maya diagrams of charge $r$,
and set $\cM_{\snz}:=
{\bigcup_{r\in \snz}}\cM_{\snz}^{(r)}$.
In particular, the element $\nz_{\leq r}\in \cM_{\snz}^{(r)}$ is called
the ground state of charge $r$. 
\vskip 1mm
\noindent
{\rm (2)} For a Maya diagram ${\bf k}$ of charge $r$, let 
${\bf k}^c:=\nz\setminus {\bf k}$ be the complement of ${\bf k}$ in $\nz$. 
We call ${\bf k}^c$ the complementary Maya diagram of charge $r$ associated
to ${\bf k}\in \cM_{\snz}^{(r)}$. We denote by $\cM_{\snz}^{(r),c}$ 
the set of all 
complementary Maya diagrams of charge $r$,
and set $\cM_{\snz}^c:=
{\bigcup_{r\in \snz}}\cM_{\snz}^{(r),c}$.
In particular, the element $(\nz_{\leq r})^c=\nz_{>r}\in \cM_{\snz}^{(r),c}$ is
called the complementary ground state of charge $r$.
\end{defn}
 
From the definition above, a Maya diagram ${\bf k}$ of charge $r$ can be 
regarded as a sequence of integers indexed by $\nz_{\leq r}$:
$${\bf k}=\{k_j~|~j\in\nz_{\leq r}\}\mbox{ such that }
k_{j-1}<k_{j}~\mbox{for }j\leq r,~ k_j=j~\mbox{for }j\ll r.$$
Similarly, the complementary Maya diagram ${\bf k}$ of charge $r$
can be regarded as a sequence of integers indexed by $\nz_{>r}$:
$${\bf k}=\{k_j~|~j\in\nz_{> r}\}\mbox{ such that }
k_j<k_{j+1}~\mbox{for }j> r,~ k_j=j~\mbox{for }j\gg r.$$
The map
$c: \cM_{\snz} \to \cM_{\snz}^c$ defined by ${\bf k}\mapsto {\bf k}^c$ is a 
bijection; the inverse of this map is also denoted by $c$.\\

Fix an integer $r\in\nz$. It is well-known that $\cM_{\snz}^{(r)}$ can
be identified with the set $\mathcal{Y}$ of all Young diagrams in the 
following way.
Let ${\bf k}=\{k_j~|~j\in\nz_{\leq r}\}\in \cM_{\nz}^{(r)}$ be a 
Maya diagram of charge $r$, and $p,q\in\nz$ integers satisfying the condition
(3.2.1). Define an increasing sequence of nonnegative integers
$0\leq m_0<m_1<\cdots<m_{p-1}\leq p+q-1$ by
$${\bf k}\cap\nz_{>r-p}
:=\{m_0+(r-p+1),~m_1+(r-p+1),~\ldots,~m_{p-1}+(r-p+1)\}.$$  
Notice that $0\leq m_0\leq m_1-1\leq\cdots\leq m_{p-1}-(p-1)\leq q$. Hence
we can define a Young diagram $Y({\bf k})\in \mathcal{Y}$ by
$$Y({\bf k}):=(\lambda_1\geq \lambda_2\geq \cdots \geq \lambda_p\geq 0),$$
where $\lambda_j:=m_{p-j}-(p-j)~\mbox{for }1\leq j\leq p$. It is easily checked 
that the definition above of $Y({\bf k})$ does not depend on the choice of 
integers $p,q$
satisfying the condition (3.2.1), and that the map ${\bf k}\mapsto Y({\bf k})$
gives a bijection from $\cM_{\snz}^{(r)}$ to $\mathcal{Y}$. We remark that,
for $\nz_{\leq r}\in \cM_{\snz}^{(r)}$, the corresponding
Young diagram $Y(\nz_{\leq r})$ is the empty set $\phi$. 

For ${\bf k}\in \cM_{\snz}^c$, we set $Y({\bf k}):=
Y({\bf k}^c)$. The map ${\bf k}\mapsto Y({\bf k})$ also gives a bijection 
from $\cM_{\snz}^{(r),c}$ to $\mathcal{Y}$.\\ 

We identify $\Xi_{\snz}$ ({\it resp}., $\Gamma_{\snz}$) with 
$\cM_{\snz}$ ({\it resp}., $\cM_{\snz}^c$); we only give
an explicit identification of $\Xi_{\snz}$ and $\cM_{\snz}$, since 
we can identify $\Gamma_{\snz}$ with $\cM_{\snz}^c$ in a similar way.

Recall that, for $\xi\in \Xi_{\snz}$, there uniquely exists
$r\in \nz$ such that $\xi=w\Lambda_r~\mbox{for }w\in W_{\snz}$. 
By the definitions, we have 
$\langle\epsilon_k,\xi\rangle_{\snz}=0\mbox{ or }1$ for each $k\in\nz$. 
For a given $\xi\in \Xi_{\snz}$, define a subset of $\nz$ by
$${\bf k}(\xi):=\{k\in\nz~|~\langle\epsilon_k,\xi\rangle_{\snz}=1\}.$$
By the construction, ${\bf k}(\xi)$ is a Maya diagram
of charge $r$, and the map $\xi\mapsto {\bf k}(\xi)$ defines a bijection
form $W_{\snz}\Lambda_r$ to $\cM_{\snz}^{(r)}$. Hence we have a bijection
$\Xi_{\snz}\overset{\sim}{\to}\cM_{\snz}$. 

Note that $\Lambda_r\in \Xi_{\snz}$ and $\Lambda_r^c\in \Gamma_{\snz}$
are identified with the (complementary) ground states:
$${\bf k}(\Lambda_r)=\nz_{\leq r}\quad\mbox{and}\quad
{\bf k}(\Lambda_r^c)=\nz_{>r}.$$

Let us recall the action of $\sigma_i,~i\in\nz,$ on $\Xi_{\snz}$ 
({\it resp}., $\Gamma_{\snz}$). Under the
identification $\Xi_{\snz}\cong\cM_{\snz}$ ({\it resp}., 
$\Gamma_{\snz}\cong\cM_{\snz}^c$), there is the induced action 
of $\sigma_i$ on $\cM_{\snz}$ ({\it resp}., $\cM_{\snz}^c$). 
It is easy to see that the explicit form of 
the action is just the transposition $(i,i+1)$ on $\nz$. 
\subsection{Definition of BZ data associated to $\nz$}
Let $I=[n+1,n+m]$ be a finite interval in $\nz.$ Set
$$\cM_{\snz}(I):=\left\{{\bf k}\in\cM_{\snz}~|~{\bf k}=\nz_{\leq n}\cup{\bf k}_I,
\mbox{ for some }{\bf k}_I\in\cM_I^{\times}\right\},$$
$$\cM_{\snz}^c(I):=\left\{{\bf k}\in\cM_{\snz}^c~|~{\bf k}={\bf k}_I\cup
\nz_{\geq n+m+2},\mbox{ for some }{\bf k}_I\in\cM_I^{\times}\right\}.$$
Define a map $\mbox{res}_I:\cM_{\snz}(I)\to \cM_I^{\times}$ by ${\bf k}\mapsto
{\bf k}_I$. It is obvious that $\mbox{res}_I$ is a bijection. 
If we set $\Omega_I({\bf k}):=
(\mbox{res}_I)^{-1}\bigl(\widetilde{I}\setminus\mbox{res}_{I}({\bf k})\bigr)$
for ${\bf k}\in \cM_{\snz}(I)$, 
then $\Omega_I({\bf k})\in \cM_{\snz}(I)$ and the map 
$\Omega_I:\cM_{\snz}(I)\to\cM_{\snz}(I)$ is also a bijection. Similarly, 
we define bijections $\mbox{res}_I^c:\cM_{\snz}^c(I)\overset{\sim}{\to}
\cM_I^{\times}$ and $\Omega_I^c:\cM_{\snz}^c(I)\overset{\sim}{\to}\cM_{\snz}^c(I)$.
The following lemma is easily verified.
\begin{lemma}\label{lemma:res}
{\rm (1)} A Maya diagram ${\bf k}$ belongs to $\cM_{\snz}(I)$ if and only if
${\bf k}^c$ belongs to $\cM_{\snz}^c(I)$.
\vskip 1mm
\noindent
{\rm (2)} For each ${\bf k}\in \cM_{\snz}(I)$, we have
$$(\mbox{\rm res}_I^c)^{-1}\bigl(\widetilde{I}\setminus
\mbox{\rm res}_I({\bf k})\bigr)={\bf k}^c,\qquad
\mbox{\rm res}_I^{-1}\bigl(\widetilde{I}\setminus
\mbox{\rm res}_I^c({\bf k}^c)\bigr)={\bf k},\eqno{(3.3.1)}$$
$$(\mbox{\rm res}_I^c)^{-1}\bigl(\mbox{\rm res}_I({\bf k})\bigr)
=\Omega_I^c({\bf k}^c),\qquad
\mbox{\rm res}_I^{-1}\bigl(\mbox{\rm res}_I^c({\bf k}^c)\bigr)
=\Omega_I({\bf k}),\eqno{(3.3.2)}$$
$$\bigl(\Omega_I({\bf k})\bigr)^c=\Omega_I^c({\bf k}^c),\qquad
\bigl(\Omega_I^c({\bf k}^c)\bigr)^c=\Omega_I({\bf k}).\eqno{(3.3.3)}$$
\end{lemma}

Let ${\bf M}=(M_{\bf k})_{{\bf k}\in \cM_{\snz}}$ be a collection of integers
indexed by $\cM_{\snz}$. For such an ${\bf M}$, we define  
${\bf M}_I:=(M_{\bf k})_{{\bf k}\in \cM_{\snz}(I)}$. By the bijection
$\mbox{res}_I:\cM_{\snz}(I)\overset{\sim}{\to} \cM_I^{\times}$, 
${\bf M}_I$ can be regarded as a collection of integers indexed by 
$\cM_I^{\times}$. Similarly, 
for ${\bf M}=(M_{\bf k})_{{\bf k}\in \cM_{\snz}^c}$, we define 
${\bf M}_I:=(M_{\bf k})_{{\bf k}\in \cM_{\snz}^c(I)}$, which is 
regarded as a collection of integers indexed by $\cM_I^{\times}$.
\begin{defn}\label{defn:BZ}
{\em (1)} A collection ${\bf M}=(M_{\bf k})_{{\bf k}\in \cM_{\snz}^c}$ of integers 
is called a complementary BZ  {\rm (}c-BZ for short{\rm )} datum associated to 
$\nz$ if it satisfies the following conditions:
\vskip 1mm
{\rm (1-a)} For each finite interval $K$ in $\nz$, 
${\bf M}_K=(M_{\bf k})_{{\bf k}\in\cM_K^{\times}}$ is an 
element of $\mathcal{BZ}_K$.
\vskip 1mm
{\rm (1-b)} For each ${\bf k}\in \cM_{\snz}^c$, there exist a finite interval 
$I$ in $\nz$ such that
\vskip 1mm
\hspace*{5mm}{\rm (1-i)} ${\bf k}\in\cM_{\snz}^c(I)$,
\vskip 1mm
\hspace*{4mm}{\rm (1-ii)} for every finite interval $J\supset I$, 
${M}_{\Omega_J^c(\bf k)}={M}_{\Omega_I^c(\bf k)}$.
\vskip 1mm
\noindent
{\rm (2)} A collection ${\bf M}=(M_{\bf k})_{{\bf k}\in \cM_{\snz}}$ of integers 
is called an $e$-BZ datum associated to $\nz$ if it satisfies the following 
conditions:
\vskip 1mm
{\rm (2-a)} For each finite interval $K$ in $\nz$, 
${\bf M}_K=(M_{\bf k})_{{\bf k}\in\cM_K^{\times}}$ is an 
element of $\mathcal{BZ}_K^{e}$.
\vskip 1mm
{\rm (2-b)} For each ${\bf k}\in \cM_{\snz}$, there exist a finite interval 
$I$ in $\nz$ such that 
\vskip 1mm
\hspace*{5mm}{\rm (2-i)} ${\bf k}\in\cM_{\snz}(I)$,
\vskip 1mm
\hspace*{4mm}{\rm (1-ii)} for every finite interval $J\supset I$, 
${M}_{\Omega_J(\bf k)}={M}_{\Omega_I(\bf k)}$.
\end{defn}
We denote by $\mathcal{BZ}_{\snz}$ ({\it resp}., $\mathcal{BZ}_{\snz}^e$) the
set of all c-BZ ({\it resp}., $e$-BZ) data associated to $\nz$.

Let us fix a c-BZ datum ${\bf M}=(M_{\bf k})_{{\bf k}\in \cM_{\snz}^c}\in 
\mathcal{BZ}_{\snz}$. For each complementary Maya diagram ${\bf k}\in 
\cM_{\snz}^c$, we denote by $\mbox{Int}^c({\bf M};{\bf k})$ the set of all 
finite intervals $I$ in $\nz$ satisfying the condition (1-b) in the definition
above. For a given $e$-BZ datum ${\bf M}=(M_{\bf k})_{{\bf k}\in \cM_{\snz}}
\in \mathcal{BZ}_{\snz}^e$, we define $\mbox{Int}^e({\bf M};{\bf k})$ similarly. 
Note that, in the latter case, ${\bf k}$ is an element of
$\cM_{\snz}.$  

\begin{rem}{\rm
Recall the identification $\cM^c_{\snz}\cong \Gamma_{\snz}$.
As mentioned before, a complementary Maya
diagram ${\bf k}\in \cM_{\snz}^c$ can be written as ${\bf k}=w\Lambda_r^c$ 
for $w\in W_{\snz}$ and $r\in \nz$. Then, the condition (1-b) in Definition 
\ref{defn:BZ} is rewritten as follows:
\vskip 1mm
(1-b)' for each $w\in W_{\snz}$ and $r\in \nz$, there exist a finite interval 
$I$ in $\nz$ such that
\vskip 1mm
\hspace*{5mm}{\rm (1-i)'} $r\in I$ and $w\in W_I$,
\vskip 1mm
\hspace*{4mm}{\rm (1-ii)'} for every finite interval $J\supset I$, 
$M_{w\Lambda_r^J}=M_{w\Lambda_r^I}$.
\vskip 1mm
\noindent
Namely, the $\mathcal{BZ}_{\snz}$ above coincides with the set of all 
BZ data of type $A_{\infty}$ in the sense of \cite{NSS}. 
}\end{rem}

For a given ${\bf M}=(M_{\bf k})_{{\bf k}\in \cM_{\snz}^c}\in 
\mathcal{BZ}_{\snz}$, introduce a new collection  
${\bf M}^*=(M^*_{\bf k})_{{\bf k}\in \cM_{\snz}}$ of integers by 
$$M^*_{\bf k}:=M_{{\bf k}^c}.$$ 
As in the finite case, the following lemma is easily verified.
\begin{lemma}
For ${\bf M}\in \mathcal{BZ}_{\snz}$, ${\bf M}^*$ is an element of 
$\mathcal{BZ}_{\snz}^e$. Moreover, the map $\ast:\mathcal{BZ}_{\snz}\to
\mathcal{BZ}_{\snz}^e$ is a bijection.
\end{lemma} 
The inverse of this bijection is also denoted by $\ast$.\\

For ${\bf M}\in \mathcal{BZ}_{\snz}$, we introduce another collection 
$\Theta({\bf M})=(\Theta(M)_{\bf k})_{{\bf k}\in \cM_{\snz}}$ of integers 
as follows. Fix ${\bf k}\in\cM_{\snz}$ and take the complement 
${\bf k}^c\in\cM_{\snz}^c$ of ${\bf k}$. 
Since ${\bf M}\in \mathcal{BZ}_{\snz}$, there exists a finite 
interval $I\in \mbox{Int}^c({\bf M};{\bf k}^c)$. Note that, for such an $I$,
we have ${\bf k}\in \cM_{\snz}(I)$ by Lemma \ref{lemma:res} (1). 
We define $\Theta(M)_{\bf k}:=M_{(\mbox{res}_I^c)^{-1}(\mbox{res}_I({\bf k}))}.$
By (3.3.2) in Lemma \ref{lemma:res}, we have 
$(\mbox{\rm res}_I^c)^{-1}\bigl(\mbox{\rm res}_I({\bf k})\bigr)
=\Omega_I^c({\bf k}^c)$. 
Therefore, by the condition (1-b) of Definition \ref{defn:BZ}, the definition 
of $\Theta({\bf M})$ does not depend on the choice of $I$.

\begin{rem}{\rm
Under the identification $\cM^c_{\snz}\cong \Gamma_{\snz}$, the bijection above
agrees with the map $\Theta$, which is introduced in \cite{NSS}. 
}\end{rem}
\subsection{Kashiwara operators on BZ data associated to $\nz$}
In \cite{NSS}, we defined the action of Kashiwara operators on 
$\mathcal{BZ}_{\snz}$. In this subsection, we reformulate some results of
\cite{NSS} under the identification $\cM^c_{\snz}\cong \Gamma_{\snz}$.\\

First, we define the action of the raising Kashiwara operators 
$\te_p$, $p\in\nz$, on $\mathcal{BZ}_{\snz}$. For ${\bf M}=
(M_{\bf k})_{{\bf k}\in \cM_{\snz}^c}\in \mathcal{BZ}_{\snz}$ and $p\in\nz$, 
define
$$\eps_p({\bf M}):=-\left(\Theta(M)_{{\bf k}(\Lambda_p)}+
\Theta(M)_{{\bf k}(\sigma_p\Lambda_p)}
-\Theta(M)_{{\bf k}(\Lambda_{p+1})}-\Theta(M)_{{\bf k}(\Lambda_{p-1})}\right).$$
By the definition of $\mathcal{BZ}_{\snz}$, $\eps_p({\bf M})$ is a nonnegative
integer. Indeed, take a finite interval $I$ from
$\mbox{Int}^c({\bf M};{\bf k}(\Lambda_p)^c)\cap 
\mbox{Int}^c({\bf M};{\bf k}(\sigma_p\Lambda_p)^c)\cap
\mbox{Int}^c({\bf M};{\bf k}(\Lambda_{p+1})^c)\cap
\mbox{Int}^c({\bf M};{\bf k}(\Lambda_{p-1})^c).$
Then, by an argument similar to the one in \cite{NSS}, we get $\eps_p({\bf M})=
\eps_p({\bf M}_I)$. Hence this is a nonnegative integer. 
 
If $\eps_p({\bf M})=0$, we set $\te_p{\bf M}={0}$. Suppose that 
$\eps_p({\bf M})>0$. Then we define $\te_p{\bf M}=
(M_{\bf k}')_{{\bf k}\in \cM_{\snz}^c}$ as follows. For ${\bf k}\in \cM_{\snz}^c$,
take a finite interval $I$ in $\nz$ such that
${\bf k}\in\cM_{\snz}^c(I)$ and 
$I\in\mbox{Int}^c({\bf M};{\bf k}(\Lambda_p)^c)\cap 
\mbox{Int}^c({\bf M};{\bf k}(\sigma_p\Lambda_p)^c)\cap
\mbox{Int}^c({\bf M};{\bf k}(\Lambda_{p+1})^c)\cap
\mbox{Int}^c({\bf M};{\bf k}(\Lambda_{p-1})^c).$
Set
$$M'_{\bf k}:=\bigl(\te_p{\bf M}_I\bigr)_{\mbox{res}_I^c({\bf k})}.$$
Here we note that $\te_p{\bf M}_I$ is defined since 
${\bf M}_I\in\mathcal{BZ}_I$.\\

Second, let us define the action of the lowering Kashiwara operators 
$\tf_p,~p\in\nz,$ on $\mathcal{BZ}_{\snz}$. For ${\bf M}=
(M_{\bf k})_{{\bf k}\in \cM_{\snz}^c}\in \mathcal{BZ}_{\snz}$ and $p\in\nz$,
we define $\tf_p{\bf M}=
(M_{\bf k}'')_{{\bf k}\in \cM_{\snz}^c}$ as follows. For ${\bf k}\in \cM_{\snz}^c$,
take a finite interval $I$ in $\nz$ such that
${\bf k}\in\cM_{\snz}^c(I)$ and 
$I\in \mbox{Int}^c({\bf M};{\bf k}(\Lambda_p)^c)\cap 
\mbox{Int}^c({\bf M};{\bf k}(\sigma_p\Lambda_p)^c).$ 
Set
$$M''_{\bf k}:=\bigl(\tf_p{\bf M}_I\bigr)_{\mbox{res}_I^c({\bf k})}.$$
\begin{prop}[\cite{NSS}]\label{prop:kas-op-NSS}
{\rm (1)} The definition above of $M'_{\bf k}$ {\rm (}{\it resp}., 
$M''_{\bf k}${\rm )} does not depend
on the choice of $I$. 
\\
{\rm (2)} For each ${\bf M}=
(M_{\bf k})_{{\bf k}\in \cM_{\snz}^c}\in \mathcal{BZ}_{\snz}$ and $p\in\nz$,
$\te_p{\bf M}$ {\rm (}{\it resp}., $\tf_p{\bf M}${\rm )} is contained in 
$\mathcal{BZ}_{\snz}\cup\{0\}$ {\rm (}{\it resp}., 
$\mathcal{BZ}_{\snz}${\rm )}.
\end{prop} 

For ${\bf M}\in \cBZ_{\snz}^e$, set $\eps_p^*({\bf M}):=\eps_p({\bf M}^*),~
p\in\nz$. Define the Kashiwara operators $\te_p^*$ and $\tf_p^*$ on 
$\cBZ_{\snz}^e$ by
$$\te_p^*{\bf M}:=\begin{cases}
\bigl(\te_p({\bf M}^*)\bigr)^*& \mbox{if }\eps_p^*({\bf M})>0,\\
0 & \mbox{if }\eps_p^*({\bf M})=0,
\end{cases}\qquad \tf_p^*{\bf M}:=\bigl(\tf_p({\bf M}^*)\bigr)^*.$$
The following corollary is easily obtained Proposition \ref{prop:kas-op-NSS}.
\begin{cor}
For each ${\bf M}\in\mathcal{BZ}_{\snz}^e$ and $p\in\nz$,
$\te_p^*{\bf M}$ {\rm (}{\it resp}., $\tf_p^*{\bf M}${\rm )} is contained in 
$\mathcal{BZ}_{\snz}^e\cup\{0\}$ {\rm (}{\it resp}., 
$\mathcal{BZ}_{\snz}^e${\rm )}.
\end{cor}
\section{Berenstein-Zelevinsky data of type $A_{l-1}^{(1)}$}
\subsection{Root datum of type $A_{l-1}^{(1)}$ }
Let us recall the notation of \cite{NSS}. Fix $l\in \nz_{\geq 3}$.
Let $\widehat{\gtg}$ be the affine Lie algebra of type $A_{l-1}^{(1)}$, 
$\widehat{\gth}$ the Cartan subalgebra of $\widehat{\gtg}$,
$\widehat{h}_i\in\widehat{\gth},~i\in\widehat{I}:=\{0,1,\cdots,l-1\}$, 
the simple coroots of $\widehat{\gtg}$, and 
$\widehat{\alpha}_i\in \widehat{\gth}^*:=
\homc(\widehat{\gth},\nc),~i\in\widehat{I}$, the simple roots of 
$\widehat{\gtg}$. Note that $\langle \widehat{h}_i,\widehat{\alpha}_j\rangle
=\widehat{a}_{ij}$ for $i,j\in\widehat{I}$. Here, $\langle\cdot,\cdot\rangle:
\widehat{\gth}\times \widehat{\gth}^*\to \nc$ is the canonical pairing, and
$\widehat{A}=(\widehat{a}_{ij})_{i,j\in\widehat{I}}$ is the Cartan matrix
of type $A_{l}^{(1)}$  with index set $\widehat{I}$; the
entries $\widehat{a}_{ij}$ are given by
$$\widehat{a}_{ij}:=\left\{\begin{array}{ll}
2 & \mbox{if }i=j,\\
-1 & \mbox{if }|i-j|=1\mbox{ or }l-1,\\
0 & \mbox{otherwise}.
\end{array}\right.$$ 
\subsection{Action of the affine Weyl group on Maya diagrams}
Consider a bijection $\tau:\nz\to\nz$ given by $\tau(j):=j+1$ for 
$j\in\nz$. It induces 
a $\nc$-linear automorphism $\tau:\gtt^*\overset{\sim}{\to}\gtt^*$
such that $\tau(\Lambda_j)=\Lambda_{j+1}$ and $\tau(\Lambda_j^c)=
\Lambda_{j+1}^c$ for all $j\in\nz$. 
Since $\alpha_j=-\Lambda_{j-1}+2\Lambda_j-\Lambda_{j+1}$,
we also have $\tau(\alpha_j)=\alpha_{j+1}$. 
Moreover, we have $\tau\circ \sigma_j=\sigma_{j+1}\circ \tau$. Hence
we have an induced group automorphism $\tau:
W_{\snz}\overset{\sim}{\to}W_{\snz}$ given by 
$\sigma_j\mapsto \tau\circ \sigma_j\circ
\tau^{-1}=\sigma_{j+1}$. Also, the restriction of $\tau:
\gtt^*\overset{\sim}{\to}\gtt^*$ to the subset
$\Xi_{\snz}$ ({\it resp}., $\Gamma_{\snz}$) gives rise to a bijection
$\tau:\Xi_{\snz}\overset{\sim}{\to}\Xi_{\snz}$ ({\it resp}., 
$\tau:\Gamma_{\snz}\overset{\sim}{\to}\Gamma_{\snz}$). 

For $i\in \widehat{I}$, define a family $S_i$ of automorphism of 
$\gtt^*$ by
$$S_i:=\{\sigma_{i+al}~|~a\in \nz\}.$$
Since $l\geq 3$,  
$\sigma_{j_1}\sigma_{j_2}=\sigma_{j_2}\sigma_{j_1}$
for all $\sigma_{j_1},\sigma_{j_2}\in S_i$, and for a fixed $\gamma\in \Xi_{\snz}$
or $\Gamma_{\snz}$, there exists a finite subset $S_i(\gamma)\subset S_i$ 
such that $\sigma_j(\gamma)=\gamma$ for all $\sigma_j\in S_i\setminus 
S_i(\gamma)$. Therefore, we can define the following infinite product 
$\widehat{\sigma}_i$ of operators acting on $\Xi_{\snz}$ and on $\Gamma_{\snz}$:
$$\widehat{\sigma}_i:=\prod_{\sigma_j\in S_i}\sigma_j.$$ 
We easily obtain the following lemma:
\begin{lemma}
{\rm (1)} For each $i\in\widehat{I}$, we have $\tau \circ \widehat{\sigma}_i=
\widehat{\sigma}_{i+1}\circ \tau$. Here we regard $i\in \widehat{I}$ as an 
element of $\nz/l\nz$. 
\vskip 1mm
\noindent
{\rm (2)} Let $\widetilde{W}_{l-1}^{(1)}$ be the group generated by 
$\widehat{\sigma}_{i},~i\in\widehat{I},$ and $\tau$. Then it is
naturally isomorphic to the extended affine Weyl group of type $A_{l-1}^{(1)}$.
Moreover, the subgroup ${W}_{l-1}^{(1)}$ of $\widetilde{W}_{l-1}^{(1)}$ 
generated by $\widehat{\sigma}_{i},~i\in\widehat{I},$ is
naturally isomorphic to the affine Weyl group of type $A_{l-1}^{(1)}$.
\end{lemma}

Let us recall the identifications $\Xi_{\snz}\cong \cM_{\snz}$ and 
$\Gamma_{\snz}\cong \cM_{\snz}^c$. Then the induced bijections 
$\tau:\cM_{\snz}^{\flat}\overset{\sim}{\to}\cM_{\snz}^{\flat}~(\flat=\phi
\mbox{ or }c)$ are given as:
\begin{align*}
\cM_{\snz}^{(r)}\ni{\bf k}=\{k_{j}~|~j\in\nz_{\leq r}\} & \mapsto 
\tau({\bf k})=\{k_{j}+1~|~j\in\nz_{\leq r}\}\in \cM_{\snz}^{(r+1)},\\
\cM_{\snz}^{(r),c}\ni{\bf k}'=\{k_{j}'~|~j\in\nz_{> r}\} & \mapsto 
\tau({\bf k}')=\{k_{j}'+1~|~j\in\nz_{> r}\}\in \cM_{\snz}^{(r+1),c}.
\end{align*}

Now, we give a quick review of the $\widetilde{W}_{l-1}^{(1)}$-action
on $\cM_{\snz}$. For details, see \cite{N} for example.
For ${\bf k}\in \cM_{\snz}$, consider an $l$-tuple of Maya diagrams
$({\bf k}^1,{\bf k}^2,\ldots,{\bf k}^{l})$ by
$${\bf k}^j:=\{k\in \nz~|~(k-1)l+j\in {\bf k}\},\quad 1\leq j\leq l.$$
It is clear that the correspondence $\cM_{\snz}\ni {\bf k}\mapsto 
({\bf k}^1,{\bf k}^2,\ldots,{\bf k}^{l})\in \bigl(\cM_{\snz}\bigr)^l$
is bijective. Moreover, if ${\bf k}$ is of charge $r$ and ${\bf k}^j$ is
of charge $r_j$, then $r=r_1+\cdots+r_l$. The action of 
$\widetilde{W}_{l-1}^{(1)}$ on $\cM_{\snz}$ can be translated as follows:
\begin{lemma}
We have
\begin{align*}
\widehat{\sigma}_0({\bf k})&=(\tau({\bf k}^l),{\bf k}^2,\ldots,{\bf k}^{l-1},
\tau^{-1}({\bf k}^1)),\\
\widehat{\sigma}_i({\bf k})&=({\bf k}^1,\ldots,{\bf k}^{i+1},{\bf k}^{i},
\ldots,{\bf k}^l),\quad 1\leq i\leq l-1,\\
\pi({\bf k})&=(\tau({\bf k}^l),{\bf k}^1,\ldots,{\bf k}^{l-1}).
\end{align*}
\end{lemma}
\begin{defn}\label{def:core}
A Maya diagram ${\bf k}\in \cM_{\snz}^{(r)}$ is called an $l$-core of charge $r$
if ${\bf k}^j$ is 
the ground state of charge $r_j$ for all $1\leq j\leq l$.
\end{defn}
\begin{prop}
A Maya diagram ${\bf k}\in \cM_{\snz}^{(r)}$ of charge $r$ belongs to the
${W}_{l-1}^{(1)}$-orbit through the ground state $\nz_{\leq r}$ 
if and only if ${\bf k}$ is an $l$-core of charge $r$.
\end{prop}

Let $Y=(\lambda_1\geq \lambda_2\geq \cdots\geq \lambda_n)\in\mathcal{Y}$ be
a Young diagram. Then $Y$ is realized as a collection of boxes arranged in
left-justified rows, with $\lambda_i$-boxes in the $i$-th row. Each box in $Y$
determines a hook, which consists of the box itself and all those boxes 
in its row to the right of the box or in its column below the box. 
The hook length of a box is 
the number of boxes in its hook. Then the following fact is well-known.
\begin{prop}
A Maya diagram ${\bf k}\in \cM_{\snz}^{(r)}$ is called an $l$-core of charge $r$
in the sense of Definition \ref{def:core} if and only if the corresponding
Young diagram $Y({\bf k})$ contains no hook whose length is a multiple of $l$. 
\end{prop}
Therefore, we obtain the following corollary.
\begin{cor}
A Maya diagram ${\bf k}\in \cM_{\snz}^{(r)}$ of charge $r$ belongs to the
${W}_{l-1}^{(1)}$-orbit through the ground state $\nz_{\leq r}$ 
if and only if the corresponding
Young diagram $Y({\bf k})$ contains no hook whose length is a multiple of $l$.
\end{cor}

\begin{rem}{\rm
(1) Since an element of $W_{l-1}^{(1)}$ is an infinite product of elements
of $W_{\snz}$, $W_{l-1}^{(1)}$ is not a subgroup of $W_{\snz}$.
\vskip 1mm
\noindent
(2) The set $\cM_{\snz}^{(r)}=W_{\snz}\bigl({\bf k}(\Lambda_r)\bigr)$ has 
infinitely many $W_{l-1}^{(1)}$-orbits.
}\end{rem}
\subsection{BZ data of type $A_{l-1}^{(1)}$}
We set $\sigma:=\tau^l$. For ${\bf M}\in \cBZ_{\snz}$, define new 
collections $\sigma({\bf M})$ and 
$\sigma^{-1}({\bf M})$ of integers indexed by $\cM_{\snz}^c$ by
$\sigma({\bf M})_{{\bf k}}:={\bf M}_{\sigma^{-1}({\bf k})}$ and 
$\sigma^{-1}({\bf M})_{{\bf k}}:={\bf M}_{\sigma({\bf k})}$ for each
${\bf k}\in \cM_{\snz}^c$, respectively. It is shown in \cite{NSS} that
both $\sigma({\bf M})$ and $\sigma^{-1}({\bf M})$ are elements of 
$\cBZ_{\snz}$. 

Similarly, for ${\bf M}\in \cBZ_{\snz}^{e}$, we can
define new collections $\sigma^{\pm}({\bf M})$, and prove that they are
elements of $\cBZ_{\snz}^{e}$.
\begin{lemma}[\cite{NSS}]
{\rm (1)} On $\cBZ_{\snz}$, we have $\Theta\circ \sigma=\sigma\circ\Theta$.
\vskip 1mm
\noindent
{\rm (2)} For ${\bf M}\in \cBZ_{\snz}$ and $p\in\nz$, 
$\eps_p(\sigma({\bf M}))=\eps_{\sigma^{-1}(p)}({\bf M})$.
\vskip 1mm
\noindent
{\rm (3)} There hold $\sigma\circ\te_p=\te_{\sigma(p)}\circ\sigma$ and 
$\sigma\circ\tf_p=\tf_{\sigma(p)}\circ\sigma$ on $\cBZ_{\snz}\cup\{{0}\}$
for all $p\in\nz$. Here it is understood that $\sigma({0})={0}$.
\end{lemma}
\begin{defn}
Set
$$\cBZ_{\snz}^{\sigma}:=\{{\bf M}\in \cBZ_{\snz}~|~
\sigma({\bf M})={\bf M}\}\quad\mbox{and}\quad
(\cBZ_{\snz}^e)^{\sigma}:=\{{\bf M}\in \cBZ_{\snz}^c~|~
\sigma({\bf M})={\bf M}\}.$$
An element ${\bf M}$ of $\cBZ_{\snz}^{\sigma}$ {\rm (}{\it resp}., 
$(\cBZ_{\snz}^e)^{\sigma}${\rm )}  
is called a $c$-BZ {\rm (}{\it resp}., 
$e$-BZ {\rm )}datum of type $A_{l-1}^{(1)}$.
\end{defn}

Now we define a crystal structure on $\cBZ_{\snz}^{\sigma}$, 
following \cite{NSS}. For ${\bf M}\in \cBZ_{\snz}^{\sigma}$ and $p\in
\widehat{I}$, we set
$$\mbox{wt}({\bf M}):=\sum_{p\in\widehat{I}}\Theta({\bf M})_{{\bf k}(\Lambda_p)}
\widehat{\alpha}_p,\quad \widehat{\eps}_p({\bf M}):=\eps_p({\bf M}),\quad
\widehat{\vphi}_p({\bf M}):=\widehat{\eps}_p({\bf M})+
\langle \widehat{h}_p,\mbox{wt}({\bf M})\rangle.$$

In order to define the action of Kashiwara operators, we need the following.
\begin{lemma}[\cite{NSS}]\label{lemma:L}
Let $q,q'\in \nz$ with $|q-q'|\geq 2$. Then, we have $\te_q\te_{p'}
=\te_{q'}\te_q$ and $\tf_q\tf_{q'}=\tf_{q'}\tf_q$, as operators from 
$\cBZ_{\snz}$ to $\cBZ_{\snz}\cup\{{0}\}$.
\end{lemma}
For ${\bf M}\in \cBZ_{\snz}^{\sigma}$ and $p\in\widehat{I}$, 
we define $\hte_p{\bf M}$ and $\htf_p{\bf M}$ as follows.
If $\widehat{\eps}_p({\bf M})=0$, we set $\hte_p{\bf M}:={0}$. If
$\widehat{\eps}_p({\bf M})>0$, then we define a new collection $\hte_p{\bf M}=
(M'_{\bf k})$ of integers indexed by $\cM_{\snz}^c$ by 
$$M'_{\bf k}:=\bigl(e_{L({\bf k},p)}{\bf M}\bigr)_{\bf k}\quad\mbox{for each }
{\bf k}\in \cM_{\snz}^c.$$
Here, $L({\bf k},p):=\{q\in p+l\nz~|~q\in {\bf k}\mbox{ and }q+1\not\in 
{\bf k}\}$ and 
$e_{L({\bf k},p)}:=\prod_{q\in L({\bf k},p)}\te_q$. By the definition, 
$L({\bf k},p)$ is a finite set such that $|q-q'|>2$ for all $q,q'\in 
L({\bf k},p)$ with $q\ne q'$. Therefore, by Lemma \ref{lemma:L}, 
$e_{L({\bf k},p)}$ is a well-defined operator on $\cBZ_{\snz}$. 

A collection $\htf_p{\bf M}=(M''_{\bf k})$ of integers indexed by $\cM_{\snz}^c$ 
is defined by 
$$M''_{\bf k}:=\bigl(f_{L({\bf k},p)}{\bf M}\bigr)_{\bf k}\quad\mbox{for each }
{\bf k}\in \cM_{\snz}^c,$$
where $f_{L({\bf k},p)}:=\prod_{q\in L({\bf k},p)}\tf_q$. By the same reasoning
as above, we see that $f_{L({\bf k},p)}$ is a well-defined operator on 
$\cBZ_{\snz}$.
\begin{prop}[\cite{NSS}]
{\rm (1)} We have $\hte_p{\bf M}\in \cBZ_{\snz}^{\sigma}\cup\{{0}\}$ 
and $\htf_p{\bf M}\in \cBZ_{\snz}^{\sigma}$.
\vskip 1mm
\noindent
{\rm (2)} The set $\cBZ_{\snz}^{\sigma}$, equipped with the maps
$\mbox{\rm wt},\widehat{\eps}_p,\widehat{\vphi}_p,\hte_p,\htf_p$, is a crystal
of type $A_{l-1}^{(1)}$.
\end{prop}

Let ${\bf O}$ be a collection of integers indexed by $\cM_{\snz}^c$ whose
${\bf k}$-component is equal to $0$ for all ${\bf k}\in \cM_{\snz}^c$. 
It is obvious that ${\bf O}\in \cBZ_{\snz}^{\sigma}$. Let 
$\cBZ_{\snz}^{\sigma}({\bf O})$ denote the connected component of the
crystal $\cBZ_{\snz}^{\sigma}$ containing ${\bf O}$. The following is
the main result of \cite{NSS}.
\begin{thm}[\cite{NSS}]\label{thm:NSS-main}
As a crystal, $\left(\cBZ_{\snz}^{\sigma}({\bf O});\mbox{\rm wt},
\widehat{\eps}_p,\widehat{\vphi}_p,\hte_p,\htf_p\right)$ is isomorphic to
$B(\infty)$ of type $A_{l-1}^{(1)}$.
\end{thm}

In a way similar to the finite case, we can define a crystal structure on 
$(\cBZ_{\snz}^{e})^{\sigma}$. By the construction, it is easy to see that 
$\ast\circ \sigma=\sigma\circ\ast$. Therefore, the restriction of 
$\ast:\cBZ_{\snz}\overset{\sim}{\to}\cBZ_{\snz}^e$ to
the subset $\cBZ_{\snz}^{\sigma}$ gives rise to a bijection
$\ast:\cBZ_{\snz}^{\sigma}\overset{\sim}{\to}(\cBZ_{\snz}^e)^{\sigma}$.
We denote by ${\bf O}^*$ the image of ${\bf O}\in \cBZ_{\snz}^{\sigma}$
under the bijection $\ast$. Then, ${\bf O}^*$ is  
a collection of integers indexed by $\cM_{\snz}$ whose
${\bf k}$-component is equal to $0$ for all ${\bf k}\in \cM_{\snz}$.

For ${\bf M}=(\mathcal{BZ}_{\snz}^e)^{\sigma}$ and $p\in\nz$, we define
$$\mbox{wt}({\bf M}):=\mbox{wt}({\bf M}^*),\quad 
\widehat{\eps}_p^*({\bf M}):=\widehat{\eps}_p({\bf M}^*),\quad
\widehat{\vphi}_p^*({\bf M}):=\widehat{\eps}_p^*({\bf M})+
\langle \widehat{h}_p,\mbox{wt}({\bf M})\rangle$$
and
$$\hte_p^*{\bf M}:=\left\{\begin{array}{ll}
(\hte_p({\bf M}^*))^* & \mbox{if }\widehat{\eps}_p^*({\bf M})>0,\\
{0} & \mbox{if }\widehat{\eps}_p^*({\bf M})=0,
\end{array}\right.\qquad
\htf_p^*:=(\htf_p({\bf M}^*))^*.$$
The following is an easy consequence of Theorem \ref{thm:NSS-main}.
\begin{cor}
{\rm (1)} The set $(\cBZ_{\snz}^{e})^{\sigma}$, equipped with the maps
$\mbox{\rm wt},\widehat{\eps}_p^*,\widehat{\vphi}_p^*,\hte_p^*,\htf_p^*$, is 
a crystal of type $A_{l-1}^{(1)}$.\\
{\rm (2)} Let
$(\cBZ_{\snz}^{e})^{\sigma}({\bf O}^*)$ denote the connected component of the
crystal $(\cBZ_{\snz}^{e})^{\sigma}$ containing ${\bf O}^*\in 
(\cBZ_{\snz}^{e})^{\sigma}$. Then, 
$\left((\cBZ_{\snz}^{e})^{\sigma}({\bf O}^*);\mbox{\rm wt},\widehat{\eps}_p^*,
\widehat{\vphi}_p^*,\hte_p^*,\htf_p^*\right)$ is isomorphic to
$B(\infty)$ of type $A_{l-1}^{(1)}$.
\end{cor}
\section{Lusztig data of infinite size and a crystal structure on them}
\subsection{Definition of Lusztig data}
\begin{defn}
{\em (1)} Let $\Delta^+_{\snz}:=\{(i,j)\in\nz\times\nz~|~i<j~\}$. 
A collection ${\bf a}=(a_{i,j})_{(i,j)\in\Delta^+_{\snz}}$ of nonnegative integers 
indexed by $\Delta^+_{\snz}$
is called a Lusztig datum associated to $\nz$ if there exist $N=N_{\bf a}>0$ 
such that
$$a_{i,j}=0\quad \mbox{for}\quad j-i\geq N_{\bf a}.\eqno{(5.1.1)}$$
We denote by ${\mathcal B}_{\snz}$ the set of all Lusztig data 
associated to $\nz$.
\vskip 1mm
\noindent
{\em (2)} For $l\in\nz_{\geq 3}$, a Lusztig datum 
${\bf a}=(a_{i,j})_{(i,j)\in\Delta^+_{\snz}}\in\cB_{\snz}$ is said to be of type 
$A_{l-1}^{(1)}$ if 
$$a_{i,j}=a_{i+l,j+l}\mbox{ for all }(i,j)\in\Delta_{\snz}^+.\eqno{(5.1.2)}$$
We denote by ${\mathcal B}_{l-1}^{(1)}$ the set of all Lusztig data of type 
$A_{l-1}^{(1)}$.
\vskip 1mm
\noindent
{\em (3)} A Lusztig datum ${\bf a}\in {\mathcal B}_{l-1}^{(1)}$ 
is said to be aperiodic if the following conditions are satisfied{\em :} 
for each $(i,j)\in\Delta_{\snz}^+$, there exists at least one $0$ in the 
$l$-tuple
of nonnegative integers
$$\{a_{i,j},a_{i+1,j+1},\ldots,a_{i+l-1,j+l-1}\}.$$
We denote by 
${\mathcal B}_{l-1}^{(1),ap}$ the set of all aperiodic Lusztig data.
\end{defn}

The set $\cB_{l-1}^{(1)}$ can be identified with the set of {\it multisegments}.
\begin{defn}
{\rm (1)} A segment of length $r$ over $\nz/l\nz$ is a sequence of 
$r$ consecutive elements in
$\nz/l\nz$:
\begin{center}
\setlength{\unitlength}{1mm}
\begin{picture}(120,10)
\put(42,7){\line(1,0){36}}\put(42,1){\line(1,0){36}}
\put(42,1){\line(0,1){6}}\put(48,1){\line(0,1){6}}
\put(54,1){\line(0,1){6}}
\put(72,1){\line(0,1){6}}\put(78,1){\line(0,1){6}}
\put(43.5,3){\small $x_1$}\put(49.5,3){\small $x_2$}
\put(60.5,3){$\cdots$}\put(73.5,3){\small $x_r$}
\put(80,1){,}
\end{picture}
\end{center}
\noindent
where $x_p=i+p-1,~1\leq p\leq r,$ for some $i\in\nz/l\nz$.
\vskip 1mm
\noindent
{\rm (2)} A multisegment over $\nz/l\nz$ is a multiset of segments
over $\nz/l\nz$. 
The set of all
multisegments over $\nz/l\nz$ is denoted by ${\bf Seg}(\nz/l\nz)$.
\end{defn}

To $(i,j)\in \Delta_{\snz}^+$, we associate the segment
of length $r=j-i$ with $x_1=i~\mbox{mod}~l\nz$. Thus we can construct 
a bijection form $\cB_{l-1}^{(1)}$ to ${\bf Seg}(\nz/l\nz)$. Note that, 
for each $(i,j)\in \Delta_{\snz}^+$, $a_{i,j}$ is just the multiplicity of the
corresponding segment.

Under the above identification $\cB_{l-1}^{(1)}\cong {\bf Seg}(\nz/l\nz)$, 
${\bf a}\in \cB_{l-1}^{(1)}$ is aperiodic if and only if the corresponding
multisegment is aperiodic in the sense of Lusztig (\cite{L2}).
\subsection{Kashiwara operators on $\cB_{\snz}$}
For ${\bf a}\in \cB_{\snz}$ and $p\in\nz$, we set
$$A_k^{(p)}({\bf a}):=\sum_{s\leq k}(a_{s,p+1}-a_{s-1,p})\quad \mbox{for }
k\leq p,$$
$$A_k^{*(p)}({\bf a}):=\sum_{t\geq k+1}(a_{p,t}-a_{p+1,t+1})\quad
\mbox{for }k\geq p.$$ 
Because of the finiteness condition (5.1.1), each of the sums above is a 
finite sum, and there exist $k_1$ and $k_2$ such that
$$A_k^{(p)}({\bf a})=0\quad\mbox{for }k\leq k_1,\qquad
A_k^{*(p)}({\bf a})=0\quad\mbox{for }k\geq k_2.$$
If we set
$$\eps_p({\bf a}):=
\mbox{max}\left\{\left.A_k^{(p)}({\bf a})~\right|~k\leq p\right\},\qquad
\eps_p^*({\bf a}):=
\mbox{max}\left\{\left.A_k^{*(p)}({\bf a})~\right|~k\geq p\right\},$$
then these are nonnegative integers. Also, set
$$\mathcal{K}(p;{\bf a}):=\left\{k~\left|~k\leq p,~\eps_p({\bf a})
=A_k^{(p)}({\bf a})\right.\right\},
\quad
\mathcal{K}^*(p;{\bf a}):=\left\{k~\left|~k\geq p,~\eps_p^*({\bf a})
=A_k^{*(p)}({\bf a})\right.\right\}.$$
Since $\mathcal{K}(p;{\bf a})$ is bounded above, we can define 
$k_f=k_f(p;{\bf a}):=\mbox{max}\{k~|~k\in \mathcal{K}(p;{\bf a})\}.$
Note that $\mathcal{K}(p;{\bf a})$ is not bounded below in general. 
However, if $\eps_p({\bf a})>0$, then it is a finite set. Therefore, 
we can define
$k_e=k_e(p;{\bf a}):=\mbox{min}\{k~|~k\in \mathcal{K}(p;{\bf a})\}$
only for ${\bf a}$ with $\eps_p({\bf a})>0$.
Similarly, we define $k_e^*=k_e^*(p;{\bf a}):=\mbox{max}\left\{k~|~k\in 
\mathcal{K}^*(p;{\bf a})\right\}$
only for ${\bf a}$ with $\eps_p^*({\bf a})>0$, and 
$k_f^*=k_f^*(p;{\bf a}):=\mbox{min}\left\{k~|~k\in 
\mathcal{K}^*(p;{\bf a})\right\}$ for arbitrary ${\bf a}$.

Now we define Kashiwara operators $\te_p$, $\tf_p$, $\te^*_p$, $\tf_p^*$,  
$p\in\nz$, on $\cB_{\snz}$ in a way similar to the finite case.\\

Let $I=[n+1,n+m]$ be a finite interval in $\nz$. Then, the set 
$\Delta_I^+=\{(i,j)~|~i,j\in\widetilde{I}\mbox{ with }i<j\}$ is naturally
regarded as a subset of $\Delta_{\snz}^+$. For a Lusztig datum 
${\bf a}=(a_{i,j})_{(i,j)\in\Delta^+_{\snz}}\in \cB_{\snz}$, we set 
${\bf a}^I:=(a_{i,j})_{(i,j)\in\Delta^+_I}$; it is an element of $\cB_I$.

Let $p\in\nz$ and ${\bf a}\in\cB_{\snz}$. Then, for each finite 
interval $I$ with $I\supset [p-N_{\bf a},p+N_{\bf a}]$, we have
$$\eps_p({\bf a})=\eps_p({\bf a}^I)\quad\mbox{and}\quad
\eps_p^*({\bf a})=\eps_p^*({\bf a}^I).\eqno{(5.2.1)}$$
Take a sufficiently large interval $I$, and assume that
$\eps_p({\bf a})=\eps_p({\bf a}^I)>0$. 
Write $\te_p{\bf a}^I=(a'_{i,j})_{(i,j)\in\Delta^+_I}\in \cB_I$, 
and define a new collection 
$\mbox{Ind}_I^{\snz}(\te_p{\bf a}^I)=(a''_{i,j})_{(i,j)\in\Delta^+_{\snz}}$
of nonnegative integers by
$$a''_{i,j}:=\begin{cases}
a'_{i,j} & \mbox{if }(i,j)\in\Delta_I^+,\\
a_{i,j} & \mbox{otherwise}.
\end{cases}$$
We also define $\mbox{Ind}_I^{\snz}(\tf_p{\bf a}^I)$, 
$\mbox{Ind}_I^{\snz}(\te_p^*{\bf a}^I)$, and 
$\mbox{Ind}_I^{\snz}(\tf_p^*{\bf a}^I)$ in a similar way.
It is obvious that 
$$\te_p{\bf a}=\mbox{Ind}_I^{\snz}(\te_p{\bf a}^I),\quad
\tf_p{\bf a}=\mbox{Ind}_I^{\snz}(\tf_p{\bf a}^I),\quad
\te_p^*{\bf a}=\mbox{Ind}_I^{\snz}(\te_p^*{\bf a}^I),\quad
\tf_p^*{\bf a}=\mbox{Ind}_I^{\snz}(\tf_p^*{\bf a}^I).\eqno{(5.2.2)}$$

\begin{rem}{\rm
Let $\cB_{\snz}^{fin}$ be the set of those Lusztig data 
${\bf a}=(a_{i,j})_{(i,j)\in \Delta_{\snz}^+}$ associated to $\nz$ for which 
$$a_{i,j}=0\quad\mbox{except for finitely many }(i,j)\in \Delta_{\snz}^+.
\eqno{(5.2.3)}$$
The (proper) subset $\cB_{\snz}^{fin}$ of $\cB_{\snz}$
should be called the set of Lusztig data of type
$A_{\infty}$. In the following, we will explain the reason. 

A Lusztig datum ${\bf a}\in\cB_{\snz}^{fin}$ can be regarded as an element 
of $\cB_I$ for a sufficiently large interval $I$, and the set
$\cB_{\snz}^{fin}\cup\{0\}$ is stable under the Kashiwara operators 
$\te_p,\tf_p,\te_p^*,\tf_p^*,~p\in\nz$. Furthermore, we can define a weight of
${\bf a}\in\cB_{\snz}^{fin}$ by
$$\mbox{wt}({\bf a}):=-\sum_{p\in\snz}
\left(\sum_{s\leq p}\sum_{t\geq p+1}a_{s,t}\right){\alpha}_p.$$
We remark that the right-hand side is a finite sum because of the finiteness
condition (5.2.3).
Hence each of 
$\left(\cB_{\snz}^{fin};\mbox{wt},\eps_p,\vphi_p,\te_p,\tf_p\right)$ and 
$\left(\cB_{\snz}^{fin};\mbox{wt},\eps_p^*,\vphi_p^*,\te_p^*,\tf_p^*\right)$ is
a crystal of type $A_{\infty}$. Moreover, they are both isomorphic to
$B(\infty)$ of type $A_{\infty}$. Indeed, it can be checked directly that 
they satisfy the conditions which uniquely characterize $B(\infty)$  
(\cite{KS}). However, since we do not use $\cB_{\snz}^{fin}$ for the remainder 
of this paper, we omit the details.

Enomoto and Kashiwara \cite{EK} gave a combinatorial description of 
$B(\infty)$ of type $A_{\infty}$ by using the PBW basis.
The crystal structure of $B(\infty)$ introduced by them agrees
with our
$\left(\cB_{\snz}^{fin};\mbox{wt},\eps_p^*,\vphi_p^*,\te_p^*,\tf_p^*\right)$. 
}\end{rem}
\subsection{Crystal structure on $\cB_{l-1}^{(1)}$ }
Let us define a crystal structure of type $A_{l-1}^{(1)}$ on $\cB_{l-1}^{(1)}$. 
The following lemma is easily shown.
\begin{lemma}\label{lemma:periodic}
Let $p\in \widehat{I}=\{0,1,\cdots,l-1\}$ and ${\bf a}\in \cB_{l-1}^{(1)}$. 
For each $r\in\nz$, we have 
$$A_k^{(p)}=A_{k+rl}^{(p+rl)},\quad A_k^{*(p)}=A_{k+rl}^{*(p+rl)},\quad
\eps_p({\bf a})=\eps_{p+rl}({\bf a}),\quad 
\eps_p^*({\bf a})=\eps_{p+rl}^*({\bf a})$$
and
$$k_e(p+rl;{\bf a})=k_e(p;{\bf a})+rl,\qquad 
k_f(p+rl;{\bf a})=k_f(p;{\bf a})+rl,$$
$$k_e^*(p+rl;{\bf a})=k_e^*(p;{\bf a})+rl,\qquad 
k_f^*(p+rl;{\bf a})=k_f^*(p;{\bf a})+rl.$$
\end{lemma}
\vskip 3mm
We now set
$$\mbox{wt}({\bf a}):=-\sum_{p\in\widehat{I}}r_p({\bf a})\widehat{\alpha}_p,
\qquad\mbox{where }r_p({\bf a}):=\sum_{s\leq p}\sum_{t\geq p+1}a_{s,t},$$
$$\widehat{\eps}_p({\bf a}):=\eps_p({\bf a}),\quad 
\widehat{\eps}_p^*({\bf a}):=\eps_p^*({\bf a}),\quad
\widehat{\vphi}_p:=\widehat{\eps}_p({\bf a})+
\langle \widehat{h}_p,\mbox{wt}({\bf a})\rangle,\quad
\widehat{\vphi}_p^*:=\widehat{\eps}_p^*({\bf a})+
\langle \widehat{h}_p,\mbox{wt}({\bf a})\rangle.$$
Note that, because of the finiteness condition (5.1.1), the right-hand side 
of the definition of $r_p({\bf a})$ is a finite sum. If we take a sufficiently 
large interval $I$, then the following is obvious from the definitions:
$$\langle \widehat{h}_{p},\mbox{wt}({\bf a})\rangle=
\langle h_p,\mbox{wt}({\bf a}^I)\rangle.\eqno{(5.3.1)}$$
Moreover, the next formulas follow immediately: 
$$\te_q\te_{q'}=\te_{q'}\te_q,\quad \tf_q\tf_{q'}=\tf_{q'}\tf_q,\quad
\te_q^*\te_{q'}^*=\te_{q'}^*\te_q^*,\quad \tf_q^*\tf_{q'}^*=\tf_{q'}^*\tf_q^*$$
for $q,q'\in\nz$ with $|q-q'|>2$. Therefore, we can define (well-defined)
operators by
$$\hte_p:=\prod_{r\in\nz}\te_{p+rl},\quad 
\htf_p:=\prod_{r\in\nz}\tf_{p+rl},\quad
\hte_p^*:=\prod_{r\in\nz}\te_{p+rl}^*,
\quad \htf_p^*:=\prod_{r\in\nz}\tf_{p+rl}^*.
$$ 
By Lemma \ref{lemma:periodic}, the image of ${\bf a}\in\cB_{l-1}^{(1)}$ under
each of the operators above belongs to $\cB_{l-1}^{(1)}\cup\{{0}\}$.
\begin{prop}
Each of $\left(\cB_{l-1}^{(1)};\mbox{\rm wt},\widehat{\eps}_p,\widehat{\vphi}_p,
\hte_p,\htf_p\right)$ and $\left(\cB_{l-1}^{(1)};\mbox{\rm wt},
\widehat{\eps}_p^*,\widehat{\vphi}_p^*,\hte_p^*,\htf_p^*\right)$ is
 a crystal of type $A_{l-1}^{(1)}$.
\end{prop}

Since one can easily check the axioms of a crystal 
(see \cite{HK}, for example) by using the definitions, 
we omit the details of the proof of this proposition.\\

\begin{defn}
A Lusztig datum ${\bf a}\in \cB_{l-1}^{(1)}$ is called a maximal element if 
$\heps_p({\bf a})=0$ for all $p\in\widehat{I}$.
\end{defn}

Let ${\bf z}=(z_1,z_2,\cdots)$ be an infinite series of nonnegative integers
such that $z_n=0$ for $n\gg 1$. Denote by $\mathcal{Z}$ the set of all such 
${\bf z}$. For ${\bf z}\in\mathcal{Z}$, we define ${\bf a}_{\bf z}=
((a_{\bf z})_{i,j})_{(i,j)\in\Delta_{\snz}^+}\in\cB_{l-1}^{(1)}$ by
$(a_{\bf z})_{i,j}:=z_{j-i}.$
For ${\bf a}\in \cB_{l-1}^{(1)}$, we introduce an infinite series 
${\bf z}({\bf a})=(z({\bf a})_1,z({\bf a})_2,\cdots)\in\mathcal{Z}$
of nonnegative integers by 
$z({\bf a})_n:=\mbox{min}\{a_{i,j}~|~j-i=n~\}$ for $n\geq 1$,
and set
$\cB_{l-1}^{(1)}({\bf z}):=\{{\bf a}\in \cB_{l-1}^{(1)}~|~{\bf z}({\bf a})
={\bf z}\}.$
Then the following are obvious from the definitions:
$$\cB_{l-1}^{(1)}=\bigsqcup_{{\bf z}\in\mathcal{Z}}\cB_{l-1}^{(1)}({\bf z})
\eqno{(5.3.2)}$$
and
$$\cB_{l-1}^{(1)}({\bf 0})=\cB_{l-1}^{(1),ap},\quad\mbox{where we set } 
{\bf 0}:=(0,0,\cdots)\in\mathcal{Z}.$$

\begin{lemma} 
{\rm (1)} Denote by $\mbox{\rm Max}(\cB_{l-1}^{(1)})$ the set of all maximal
elements in $\cB_{l-1}^{(1)}$. Then we have $\mbox{\rm Max}
(\cB_{l-1}^{(1)})=\{{\bf a}_{\bf z}~|~{\bf z}\in\mathcal{Z}\}.$ 
\vskip 1mm
\noindent
{\rm (2)} A Lusztig datum ${\bf a}\in\cB_{l-1}^{(1)}$ is a maximal element if 
and only if $\heps_p^*({\bf a})=0$ for all $p\in\widehat{I}$.
\end{lemma} 
\begin{proof} By the definitions, we have ${\bf a}_{\bf z}\in 
\mbox{Max}(\cB_{l-1}^{(1)})$. Conversely, let ${\bf a}\in \cB_{l-1}^{(1)}$ 
be a maximal element. By the decomposition (5.3.2), there exists a 
unique ${\bf z}\in\mathcal{Z}$ such that ${\bf a}\in \cB_{l-1}^{(1)}({\bf z})$. 
What we need to prove is that ${\bf a}={\bf a}_{\bf z}$.
Introduce a new 
Lusztig datum ${\bf a}(0)=(a(0)_{i,j})\in \cB_{l-1}^{(1)}({\bf 0})$ by
$a(0)_{i,j}:=a_{i,j}-z_{j-i}$ for each $(i,j)\in\Delta_{\snz}^+.$
By the definitions, we have $\heps_p({\bf a})=\heps_p({\bf a}(0))$ for all
$p\in\widehat{I}$. Therefore, we may assume that ${\bf z}={\bf 0}$. 

Let ${\bf a}\in \cB_{l-1}^{(1)}({\bf 0})=\cB_{l-1}^{(1),ap}$. 
We will prove that
{\it if ${\bf a}\ne {\bf a}_{\bf 0}$, then there exists $p\in\widehat{I}$ 
such that $\widehat{\eps}_p({\bf a})>0$}. Set $B_k({\bf a})=
\mbox{max}\{a_{k,1},a_{k+1,2},\ldots,a_{k+l-1,l}\}$ for $k\leq 0$. 
Since $B_k({\bf a})=0$ for all sufficiently 
small $k$, there exists the minimum $k_0\leq 0$ for which $B_{k_0}({\bf a})>0$. 
Because of the aperiodicity condition for ${\bf a}$, 
there exists $p_0\in \widehat{I}$ 
such that $a_{k_0+p_0-1,p_0}=0$ and $a_{k_0+p_0,p_0+1}>0$. Here 
the indices are regarded as elements of $\nz/l\nz$. Hence we have 
\begin{align*}
A_{k_0+p_0}^{(p_0)}({\bf a})&=
\sum_{s\leq k_0+p_0}(a_{s,p+1}-a_{s-1,p})=a_{k_0+p_0,p_0+1}>0.
\end{align*}
Therefore, we deduce that $\widehat{\eps}_p(\bf a)=\eps_p({\bf a})>0$.
This proves part (1).

Let $\mbox{Max}^*(\cB_{l-1}^{(1)})$ be the set of those elements ${\bf a}$ in
$\cB_{l-1}^{(1)}$ such that $\heps_p^*({\bf a})=0$ for all $\widehat{I}$. 
By arguing as in the proof of part (1), we obtain $\mbox{Max}^*
(\cB_{l-1}^{(1)})=\{{\bf a}_{\bf z}~|~{\bf z}\in\mathcal{Z}\}.$ 
This proves part (2).
\end{proof}
The following corollary is an easy consequence of the lemma above.
\begin{cor}
For each ${\bf z}=(z_n)\in\mathcal{Z}$, $\cB_{l-1}^{(1)}({\bf z})$ contains a 
unique maximal element ${\bf a}_{\bf z}$ for which
$\mbox{\rm wt}({\bf a}_{\bf z})=-m({\bf z})\delta.$
Here, $m({\bf z})=\sum_{n\geq 1}nz_n$, and  
$\delta:=\sum_{p\in\widehat{I}}\widehat{\alpha}_p$ is the null root.
\end{cor}

Let $\widehat{P}$ denote the weight lattice of type $A_{l-1}^{(1)}$, and 
$\widehat{\lambda}\in \widehat{P}$. 
Consider the set $T_{\widehat{\lambda}}=\{t_{\widehat{\lambda}}\}$, 
and introduce two crystal structures 
$(T_{\widehat{\lambda}};\mbox{wt},\heps_p,\hvphi_p,\hte_p,\htf_p)$ and 
$(T_{\widehat{\lambda}};\mbox{wt},\heps_p^*,\hvphi_p^*,\hte_p^*,\htf_p^*)$
of type $A_{l-1}^{(1)}$ as follows: $\mbox{wt}(t_{\widehat{\lambda}})
={\widehat{\lambda}}$, and for each
$p\in\widehat{I}$, $\heps_p(t_{\widehat{\lambda}})=\hvphi_p(t_{\widehat{\lambda}})
=\heps_p^*(t_{\widehat{\lambda}})=\hvphi_p^*(t_{\widehat{\lambda}})=-\infty$, 
$\hte_pt_{\widehat{\lambda}}=\htf_pt_{\widehat{\lambda}}=
\hte_p^*t_{\widehat{\lambda}}=\htf_p^*t_{\widehat{\lambda}}=0$.   

Define a map $\mathcal{T}:\cB_{l-1}^{(1)}({\bf z})\to 
\cB_{l-1}^{(1)}({\bf 0})\otimes T_{-m({\bf z})\delta}$ by
$${\bf a}\quad\mapsto\quad {\bf a}(0)\otimes t_{-m({\bf z})\delta}.$$
\begin{thm}[\cite{LTV}]\label{thm:LTV}
{\rm (1)} For an arbitrary ${\bf z}\in\mathcal{Z}$, each of 
$\left(\cB_{l-1}^{(1)}({\bf z});\mbox{\rm wt},\widehat{\eps}_p,\widehat{\vphi}_p,
\hte_p,\htf_p\right)$ and 
$\left(\cB_{l-1}^{(1)}({\bf z});\mbox{\rm wt},\widehat{\eps}_p^*,
\widehat{\vphi}_p^*,\hte_p^*,\htf_p^*\right)$ is a crystal of type 
$A_{l-1}^{(1)}$.
\vskip 1mm
\noindent
{\rm (2)} Each of $\left(\cB_{l-1}^{(1)}({\bf 0});
\mbox{\rm wt},\widehat{\eps}_p,\widehat{\vphi}_p,
\hte_p,\htf_p\right)$ and 
$\left(\cB_{l-1}^{(1)}({\bf 0});\mbox{\rm wt},\widehat{\eps}_p^*,
\widehat{\vphi}_p^*,\hte_p^*,\htf_p^*\right)$ is isomorphic to $B(\infty)$ 
of type $A_{l-1}^{(1)}$.
\vskip 1mm
\noindent
{\rm (3)} The map $\mathcal{T}:\cB_{l-1}^{(1)}({\bf z})\to 
\cB_{l-1}^{(1)}({\bf 0})\otimes T_{-m({\bf z})\delta}$ is an isomorphism of 
crystals with respect to each of the crystal structures above.
\vskip 1mm
\noindent
{\rm (4)} The decomposition $(5.3.2)$ gives a decomposition of $\cB_{l-1}^{(1)}$
into its connected components. More specifically, we have 
$$\cB_{l-1}^{(1)}\cong \bigoplus_{m\in\snz_{\geq 0}}\bigl(B(\infty)\otimes 
T_{-m\delta}\bigr)^{\oplus p(m)},$$
where $p(m)$ denotes the number of partitions of $m\geq 0$.
\end{thm}

\begin{rem}{\rm 
The crystal structure on $\cB_{l-1}^{(1)}$ (or equivalently,
on ${\bf Seg}(\nz/l\nz)$) introduced in \cite{LTV} is 
$\left(\cB_{l-1}^{(1)};\mbox{\rm wt},\widehat{\eps}_p^*,
\widehat{\vphi}_p^*,\hte_p^*,\htf_p^*\right)$. Therefore, 
strictly speaking, they proved 
the statement above only for this crystal structure. 
However, by a similar method, one can 
prove the statement for the other crystal structure. Hence 
we omit the details.
}\end{rem}
\section{BZ data arising from Lusztig data of Infinite size }
\subsection{Results on BZ data associated to finite intervals}
In this subsection, we prove some results on 
BZ data associated to finite intervals, which we need later.
 
Throughout this subsection, ${\bf k}$ denotes a Maya diagram associated to
a finite interval. More specifically, let $I=[n+1,n+m]$ be a finite interval 
and ${\bf k}=(k_{n+1}<\cdots< k_{n+u})\in \cM_{I}^{\times}$. Recall that 
$k_j\in \widetilde{I}=[n+1,n+m+1]$ for all
$n+1\leq j\leq n+u$. 

Let $d\in \nz_{>0}$, and define
$$I':=[n-d+1,n+m],\quad {\bf k}':=[n-d+1,n]\cup{\bf k}
\in\cM_{I'}^{\times},$$
$$I'':=[n+1,n+m+d],\quad {\bf k}'':={\bf k}\cup 
[n+m+2,n+m+d+1]\in\cM_{I''}^{\times}.$$
We use the notation above throughout this subsection.
\begin{lemma}\label{lemma:fin1} 
Let ${\bf a}\in\cB_{\snz}$, and suppose that $u\geq N_{\bf a}$.
\vskip 1mm
\noindent
{\rm (1)} If $k_j=j$ for each $n+1\leq j\leq n+N_{\bf a}$, then we have
$M_{{\bf k}}({\bf a}^I)=M_{{\bf k}'}({\bf a}^{I'})$ for all $d\in\nz_{>0}$.
\vskip 1mm
\noindent 
{\rm (2)} If $k_j=m-u+j$ for each $n+u-N_{\bf a}+1\leq j\leq n+u$, then
we have $M_{\bf k}({\bf a}^I)=M_{{\bf k}''}({\bf a}^{I''})$ for all $d\in\nz_{>0}$.
\end{lemma}
\begin{proof}
We only give a proof of part (1), since part (2) is shown in a similar way.
We write 
$${\bf k}'=(k_{n-d+1}<\cdots<k_n<k_{n+1}<\cdots< k_{n+u}).$$
Then we have
$$k_j=j\quad \mbox{for each }n-d+1\leq j\leq n+N_{\bf a}.\eqno{(6.1.1)}$$
The explicit form of 
$M_{{\bf k}}({\bf a}^{I})$ is given by
$$
M_{{\bf k}}({\bf a}^{I})=
-\sum_{j=n+1}^{n+u}\sum_{i=n+1}^{k_j-1}a_{i,k_j}+
\mbox{min}\left\{\left.\sum_{n+1\leq p<q\leq n+u}a_{c_{p,q},c_{p,q}+(q-p)}
~\right|~\begin{array}{c}C=(c_{p,q})\mbox{ is }\\ 
\mbox{a ${\bf k}$-tableau}\end{array}\right\}.
$$
Let us compute the right-hand side. By (5.1.1) and (6.1.1), we have
\begin{align*}
\sum_{j=n+1}^{n+u}\sum_{i=n+1}^{k_j-1}a_{i,k_j}&=
\sum_{j=n+1}^{n+N_{\bf a}}\sum_{i=n+1}^{k_j-1}a_{i,k_j}+
\sum_{j=n+N_{\bf a}+1}^{n+u}\sum_{i=n+1}^{k_j-1}a_{i,k_j}\\
&=\sum_{n+1\leq i<j\leq n+N_{\bf a}}a_{i,j}+
\sum_{j=n+N_{\bf a}+1}^{n+u}\sum_{i=k_j-N_{\bf a}+1}^{k_j-1}a_{i,k_j}. 
\end{align*}
Now, let $C=(c_{p,q})_{n+1\leq p<q\leq n+u}$ be a 
${\bf k}^{I}$-tableau.
Because of (6.1.1) and the monotonicity condition for $C$ , 
we have $c_{p,q}=p$ for $n+1\leq p<q\leq n+N_{\bf a}$. 
Therefore, by (5.1.1), we compute:
\begin{align*}
&\sum_{n+1\leq p<q\leq n+u}a_{c_{p,q},c_{p,q}+(q-p)}\\
&\qquad\qquad
=\sum_{n+1\leq p<q\leq n+N_{\bf a}}a_{c_{p,q},c_{p,q}+(q-p)}+
\sum_{q=n+N_{\bf a}+1}^{n+u}\sum_{p=n+1}^{q-1}a_{c_{p,q},c_{p,q}+(q-p)}\\
&\qquad\qquad
=\sum_{n+1\leq p<q\leq n+N_{\bf a}}a_{p,q}+
\sum_{q=n+N_{\bf a}+1}^{n+u}\sum_{p=q-N_{\bf a}}^{q-1}a_{c_{p,q},c_{p,q}+(q-p)}.
\end{align*}
Consequently, we deduce that
\begin{align*}
M_{{\bf k}}({\bf a}^{I})&=
-\sum_{j=n+N_{\bf a}+1}^{n+u}\sum_{i=k_j-N_{\bf a}+1}^{k_j-1}a_{i,k_j}\\
&\qquad
+\mbox{min}\left\{\left.\sum_{q=n+N_{\bf a}+1}^{n+u}\sum_{p=q-N_{\bf a}}^{q-1}
a_{c_{p,q},c_{p,q}+(q-p)}
~\right|~\begin{array}{c}C=(c_{p,q})\mbox{ is }\\ 
\mbox{a ${\bf k}$-tableau}\end{array}\right\}.
\end{align*}
By a similar computation, we deduce that
\begin{align*}
M_{{\bf k}'}({\bf a}^{I'})&=
-\sum_{j=n+N_{\bf a}+1}^{n+u}\sum_{i=k_j-N_{\bf a}+1}^{k_j-1}a_{i,k_j}\\
&\qquad
+\mbox{min}\left\{\left.\sum_{q=n+N_{\bf a}+1}^{n+u}\sum_{p=q-N_{\bf a}}^{q-1}
a_{c'_{p,q},c'_{p,q}+(q-p)}
~\right|~\begin{array}{c}C'=(c'_{p,q})\mbox{ is }\\ 
\mbox{a ${\bf k}'$-tableau}\end{array}\right\}.
\end{align*}
Since it is obvious that
\begin{align*}
& \mbox{min}\left\{\left.\sum_{q=n+N_{\bf a}+1}^{n+u}\sum_{p=q-N_{\bf a}}^{q-1}
a_{c_{p,q},c_{p,q}+(q-p)}
~\right|~\begin{array}{c}C=(c_{p,q})\mbox{ is }\\ 
\mbox{a ${\bf k}$-tableau}\end{array}\right\}\\
&\qquad\qquad=
\mbox{min}\left\{\left.\sum_{q=n+N_{\bf a}+1}^{n+u}\sum_{p=q-N_{\bf a}}^{q-1}
a_{c'_{p,q},c'_{p,q}+(q-p)}
~\right|~\begin{array}{c}C'=(c'_{p,q})\mbox{ is }\\ 
\mbox{a ${\bf k}'$-tableau}\end{array}\right\},
\end{align*}
we obtain the desired equality.
\end{proof}

Note that ${\bf k}\in \cM_I^{\times}$ can be naturally regarded as an element
of $\cM_{I'}^{\times}$ and of $\cM_{I''}^{\times}$. The next lemma 
follows easily from the definitions and (5.1.1).
\begin{lemma}\label{lemma:fin2} 
Let ${\bf a}\in\cB_{\snz}$.
\vskip 1mm
\noindent
{\rm (1)} If $k_{n+1}-n\geq N_{\bf a}$, then $M_{{\bf k}}({\bf a}^I)
=M_{{\bf k}}({\bf a}^{I'})$ for all $d\in\nz_{>0}$.
\vskip 1mm
\noindent
{\rm (2)} For all $d\in\nz_{>0}$, we have $M_{{\bf k}}({\bf a}^I)
=M_{{\bf k}}({\bf a}^{I''})$.
\end{lemma}

\begin{lemma}\label{lemma:fin3}
Let ${\bf a}\in \cB_{l-1}^{(1),ap}$. Assume that the following condition is
satisfied: 
$$k_{n+1}-n\geq (2l+1)N_{\bf a}.\eqno{(6.1.2)}$$
Then, we have
$M_{{\bf k}}({\bf a}^I)=M_{{\bf k}'}({\bf a}^{I'})$ for all $d\in\nz_{>0}$.
\end{lemma}
\begin{proof}
We use the same notation as in the proof of Lemma \ref{lemma:fin1}: 
$${\bf k}'=(k_{n-d+1}<\cdots<k_n<k_{n+1}<\cdots< k_{n+u}).$$
In this case, we have
$$k_j=j\quad \mbox{for each }n-d+1\leq j\leq n.\eqno{(6.1.3)}$$
\vskip 3mm
Since $k_{n+1}-n\geq (2l+1)N_{\bf a}>N_{\bf a}$, we obtain
\begin{align*}
M_{{\bf k}}({\bf a}^{I})&=
-\sum_{j=n+1}^{n+u}\sum_{i=k_j-N_{\bf a}+1}^{k_j-1}a_{i,k_j}\\
&\qquad
+\mbox{min}\left\{\left.
\sum_{q=n+2}^{n+u}\sum_{p=n+1}^{q-1}a_{c_{p,q},c_{p,q}+(q-p)}
~\right|~\begin{array}{c}C=(c_{p,q})\mbox{ is }\\ 
\mbox{a ${\bf k}$-tableau}\end{array}\right\}.
\end{align*}
By (6.1.3) and an argument similar to the one in the proof of 
Lemma \ref{lemma:fin1}, we deduce that
\begin{align*}
M_{{\bf k}'}({\bf a}^{I'})&=
-\sum_{j=n+1}^{n+u}\sum_{i=k_j-N_{\bf a}+1}^{k_j-1}a_{i,k_j}\\
&\qquad
+\mbox{min}\left\{\left.\sum_{q=n+1}^{n+u}\sum_{p=n-d+1}^{q-1}
a_{c_{p,q}',c_{p,q}'+(q-p)}
~\right|~\begin{array}{c}C'=(c'_{p,q})\mbox{ is }\\ 
\mbox{a ${\bf k}'$-tableau}\end{array}\right\}.
\end{align*}
Let us write
$$S:=\mbox{min}\left\{\left.
\sum_{q=n+2}^{n+u}\sum_{p=n+1}^{q-1}a_{c_{p,q},c_{p,q}+(q-p)}
~\right|~\begin{array}{c}C=(c_{p,q})\mbox{ is }\\ 
\mbox{a ${\bf k}$-tableau}\end{array}\right\},$$
$$S':=\mbox{min}\left\{\left.\sum_{q=n+1}^{n+u}\sum_{p=n-d+1}^{q-1}
a_{c_{p,q}',c_{p,q}'+(q-p)}
~\right|~\begin{array}{c}C'=(c'_{p,q})\mbox{ is }\\ 
\mbox{a ${\bf k}'$-tableau}\end{array}\right\}.$$
By the computation above, it suffices to show that $S=S'$ under the
assumption (6.1.2). Observe that $S'$ can be rewritten as
$$S'=\mbox{min}\left\{\left.
\sum_{q=n+2}^{n+u}\sum_{p=n+1}^{q-1}a_{c_{p,q}',c_{p,q}'+(q-p)}
+\sum_{q=n+1}^{n+u}\sum_{p=n-d+1}^{n}
a_{c_{p,q}',c_{p,q}'+(q-p)}
~\right|~\begin{array}{c}C'=(c'_{p,q})\mbox{ is }\\ 
\mbox{a ${\bf k}'$-tableau}\end{array}\right\}.$$
For each ${\bf k'}$-tableau $C'=(c'_{p,q})_{n-d+1\leq p\leq q\leq n+u}$, 
the subtableau $(c'_{p,q})_{n+1\leq p\leq q\leq n+u}$ of $C'$ is always a 
${\bf k}$-tableau. Therefore, what we need to show is the following: 
if the assumption (6.1.2) holds, then there exists a ${\bf k}'$-tableau 
$C'=(c'_{p,q})_{n-d+1\leq p\leq q\leq n+u}$ such that 
$$\displaystyle{\sum_{q=n+2}^{n+u}\sum_{p=n+1}^{q-1}a_{c_{p,q}',c_{p,q}'+(q-p)}
=S}\quad\mbox{and}\quad
\displaystyle{\sum_{q=n+1}^{n+u}\sum_{p=n-d+1}^{n}a_{c_{p,q}',c_{p,q}'+(q-p)}
=0.}\eqno{(6.1.4)}$$
In the following, we will construct such a $C'$.\\

For $1\leq v\leq N_{\bf a}$, consider the $l$-tuple of nonnegative 
integers:
$$\mathcal{A}_v:=\left\{\right.a_{i,i+N_{\bf a}+1-v}~\left|~n+v\leq i\leq n+v+l-1
\right\}.$$
Note that, by assumption (6.1.2), the following inequalities hold:
\begin{align*}
&n-d+1\leq i<i+N_{\bf a}+1-v\leq n+m\\
&\qquad\qquad\mbox{for all }1\leq v\leq N_{\bf a}\mbox{ and }
n+v\leq i\leq n+v+l-1.
\end{align*}
By the aperiodicity condition for ${\bf a}$, there exists at least one $0$ 
in $\mathcal{A}_v$. Set
$$b_{v}:=\mbox{min}\left\{~i~\left|~
a_{i,i+N_{\bf a}+1-v}\in\mathcal{A}_v~\mbox{and}~a_{i,i+N_{\bf a}+1-v}=0
\right\}\right..$$
Then it follows that
$$n+v\leq b_v\leq n+v+l-1.\eqno{(6.1.5)}$$

Set $N_{0}:=\mbox{max}\{N_{\bf a},u,d\}$, and introduce an 
$(N_{0}\times N_{0})$-matrix $X=(x_{s,t})$, where the indices run over
$n-N_0+1\leq s \leq n$ and $n+1\leq t\leq n+N_0$, by
$$x_{s,t}:=\begin{cases}
s & \mbox{if } n-N_0+1\leq s\leq n-N_{\bf a},\\
b_1+2(s-n+N_{\bf a}-1)l & \mbox{if }
n-N_{\bf a}+1\leq s\leq n,~s+N_{\bf a}\leq t,\\
b_{s-t+N_{\bf a}+1}+(s+t-2n+N_{\bf a}-2)l & 
\mbox{if }n-N_{\bf a}+1\leq s\leq n,~s+N_{\bf a}>t.
\end{cases}$$
\vskip 3mm
\noindent
{\bf Claim}. (1) {\it The entries of $X$ satisfy the {\rm (}usual{\rm )}
monotonicity conditions. Namely, for all $s$ and $t$, 
$x_{s,t}\leq x_{s,t+1}$ and $x_{s,t}<x_{s+1,t}.$}
\vskip 1mm
\noindent
(2) {\it For all $n-N_0+1\leq s \leq n$, we have $s\leq x_{s,n+1}$.}
\vskip 3mm
\noindent
{\it Proof of the claim.}
If $n-N_0+1\leq s\leq n-N_{\bf a}$, or if $n-N_{\bf a}+1\leq s\leq n$ and 
$s+N_{\bf a}\leq t$, then the monotonicity conditions are obviously
satisfied. Suppose that $n-N_{\bf a}+1\leq s\leq n$ and 
$s+N_{\bf a}>t$.
By (6.1.5), we have $b_{s-t+N_{\bf a}}\geq n+s-t+N_{\bf a}$ and 
$b_{s-t+N_{\bf a}+1}\leq n+s-t+N_{\bf a}+l$. Therefore, we see that
\begin{align*}
x_{s,t+1}-x_{s,t}&=b_{s-t+N_{\bf a}}-b_{s-t+{\bf a}+1}+l\\
&\geq (n+s-t+N_{\bf a})-(n+s-t+N_{\bf a}+l)+l\\
&\geq 0.
\end{align*}
By a similar computation, we obtain $x_{s,t}<x_{s+1,t}$. This proves
part (1).

We prove part (2). If $n-N_0+1\leq s\leq n-N_{\bf a}$, then the assertion 
is obvious.
Let $n-N_{\bf a}+1\leq s\leq n$. By (6.1.5), we see that
\begin{align*}
x_{s,n+1}-s&=b_{s-(n+1)+N_{\bf a}+1}+(s+(n+1)-2n+N_{\bf a}-2)l-s\\
&\geq n+(s-n+N_{\bf a})+(s-n+N_{\bf a}-1)l-s\\
&=N_{\bf a}+(s-n+N_{\bf a}-1)l\\
&\geq 0.
\end{align*}
This proves part (2).
\hfill$\square$
\vskip 3mm

Now let us construct a ${\bf k}'$-tableau 
$C'=(c'_{p,q})_{n-d+1\leq p\leq q\leq n+u}$ satisfying (6.1.4).\\

First, we remark that $c'_{p,q}$ must be equal to $p$ for 
$n-d+1\leq p\leq q\leq n$ by (6.1.3). 
Second, for $n-d-1\leq p\leq n$ and $n+1\leq q\leq n+u$, we set
$c'_{p,q}:=x_{p,q}$. Here, $x_{p,q}$ is the $(p,q)$-th entry of the matrix 
$X$ constructed above. By the claim above, the monotonicity conditions 
are satisfied for all $n-d+1\leq p\leq n$ and $p\leq q\leq n+u$.
Finally, define $c'_{p,q}$ for $n+1\leq p\leq q\leq n+u$ as follows. 
Let $C=(c_{p,q})_{n+1\leq p\leq q\leq n+u}$ be a ${\bf k}$-tableau such that
${\sum_{q=n+2}^{n+u}\sum_{p=n+1}^{q-1}a_{c_{p,q},c_{p,q}+(q-p)}=S}$. We 
set $c'_{p,q}=c_{p,q}$ for $n+1\leq p\leq q\leq n+u$. 
The upper triangular matrix $C'=(c'_{p,q})_{n-d+1\leq p\leq q\leq n+u}$ thus 
obtained
is a ${\bf k}'$-tableau. Indeed, by the monotonicity condition for $X$, the 
largest matrix entry is $x_{n,n+N_0}=b_1+2(N_{\bf a}-1)l$. By (6.1.2) 
and (6.1.5), we see that
\begin{align*}
k_{n+1}-x_{n,n+N_0}&\geq k_{n+1}-(n+l)-2(N_{\bf a}-1)l\\
&=(k_{n+1}-n)-(2l+1)N_{\bf a}+l+N_{\bf a}\\
&>0.
\end{align*}
Hence it follows that $c'_{n,q}\leq x_{n,n+N_0}<k_{n+1}\leq c'_{n+1,q}$ for all 
$n+1\leq q\leq n+u$. Therefore, the monotonicity conditions are satisfied 
for all $n-d-1\leq p\leq q\leq n+u$.

Moreover, by the construction, the condition (6.1.4) is automatically 
satisfied. Thus, the lemma is proved.
\end{proof}
\subsection{Collections of integers arising from Lusztig data associated to
$\nz$}
Let us return to the infinite case.
\begin{cor}\label{cor:Lres}
Let ${\bf a}\in\cB_{\snz}$ and ${\bf k}\in \cM_{\snz}$. Then, there exists a
finite interval $I_0$ such that 
{\rm (i)} ${\bf k}\in \cM_{\snz}(I_0)$, and 
{\rm (ii)} for every finite interval $J$ containing $I_0$, 
$M_{\mbox{\rm res}_J({\bf k})}({\bf a}^J)
=M_{\mbox{\rm res}_{I_0}({\bf k})}({\bf a}^{I_0})$.
\end{cor}
\begin{proof}
Take a finite interval $I=[n+1,n+m]$ such that 
${\bf k}\in \cM_{\snz}(I)$. Note that
${\bf k}\in \cM_{\snz}(J)$ for all $J\supset I$. Set $I_0:=[n-N_{\bf a}+1,n+m]$.
Then, by Lemma \ref{lemma:fin1} (1), $I_0$ has the desired properties.
\end{proof}
For given ${\bf a}\in \cB_{\snz}$ and ${\bf k}\in\cM_{\snz}$, take a
finite interval $I_0$ of Corollary \ref{cor:Lres} and define 
$M_{\bf k}({\bf a}):=M_{\mbox{res}_{I_0}({\bf k})}({\bf a}^{I_0})$. 
Note that this definition does not depend on the choice of $I_0$ by the 
corollary above.
Set ${\bf M}({\bf a}):=(M_{\bf k}({\bf a}))_{{\bf k}\in\cM_{\snz}}$, and 
define a map $\Phi_{\snz}$ from $\cB_{\snz}$ to
the set of all collections of integers indexed by $\cM_{\snz}$, by the 
assignment ${\bf a}\mapsto {\bf M}({\bf a})$.

We remark that $\Phi_{\snz}$ is not injective. Indeed, take
${\bf 0}=(0,0,\ldots)\in \mathcal{Z}$ and consider ${\bf a}_{\bf 0}\in
\cB_{\snz}$. Then we have $\Phi_{\snz}({\bf a}_{\bf 0})={\bf O}^*
\in \cBZ_{\snz}^e$.
Also, take ${\bf 1}=(1,0,0,\ldots)\in \mathcal{Z}$ and 
consider ${\bf a}_{\bf 1}\in\cB_{\snz}$; namely, ${\bf a}_{\bf 1}=(a_{i,j})$
is given by $a_{i,j}=\delta_{j-i,1}$ for each $(i,j)\in \Delta_{\snz}^+$.
We may take $N_{{\bf a}_{\bf 1}}=2$. 
Let ${\bf k}=\{k_j~|~j\in\nz_{\leq r}\}\in\cM_{\snz}^{(r)}$ and 
$I_0=[n_0+1,n_0+m_0]$. Then we compute:
\begin{align*}
M_{\bf k}({\bf a}_{\bf 1})&=
M_{res_{I_0}({\bf k}^{I_0})}\left(({{\bf a}_{\bf 1}})^{I_0}\right)\\
&=-\sum_{j=n_0+N_{\bf a}+1}^{r}\sum_{i=k_j-1}^{k_j-1}\delta_{k_j-i,1}
+\mbox{min}\left\{\left.\sum_{q=n_0+N_{\bf a}+1}^{r}\sum_{p=q-2}^{q-1}
\delta_{q-p,1}
~\right|~\begin{array}{c}C=(c_{p,q})\mbox{ is }\\ 
\mbox{a ${\bf k}^{I_0}$-tableau}\end{array}\right\}\\
&=-(r-n_0-N_{\bf a})+(r-n_0-N_{\bf a})\\
&=0.
\end{align*}
Therefore,  we conclude that  $\Phi_{\snz}({\bf a}_{\bf 1})={\bf O}^*
=\Phi_{\snz}({\bf a}_{\bf 0})$. 

By a similar computation, we have the following.
\begin{lemma} 
Let ${\bf a}_{\bf z}$ be a maximal  
element in $\cB_{l-1}^{(1)}({\bf z})$. Then we have 
$\Phi_{\snz}({\bf a}_{\bf z})={\bf O}^*.$
\end{lemma}
\subsection{BZ data arising from aperiodic Lusztig data}
The aim of this subsection is to prove the following proposition.
\begin{prop}\label{prop:e-BZ}
For each ${\bf a}\in \cB_{l-1}^{(1),ap}=\cB_{l-1}^{(1)}({\bf 0})$, the collection
${\bf M}({\bf a})=(M_{\bf k}({\bf a}))_{{\bf k}\in\cM_{\snz}}$ is an $e$-BZ datum
in the sense of Definition \ref{defn:BZ}.
In other words, the restriction of $\Phi_{\snz}$ to $\cB_{l-1}^{(1),ap}$
defines a map $\Phi_{\snz}:\cB_{l-1}^{(1),ap}\to \cBZ_{\snz}^{e}$.
\end{prop}

In order to prove the proposition, we need the next lemma.
\begin{lemma}\label{lemma:BZ-3}
Let ${\bf a}\in \cB_{l-1}^{(1),ap}$ and ${\bf k}\in\cM_{\snz}$. Then, there
exists a finite interval $I$ in $\nz$ such that, for every finite interval
$J\supset I$, $M_{\Omega_J({\bf k})}({\bf a})=M_{\Omega_I({\bf k})}({\bf a})$.
\end{lemma}
\begin{proof} We may assume that ${\bf k}$ is a Maya diagram of charge $r$, 
and write ${\bf k}=\{k_j~|~j\in\nz_{\leq r}\}$. Define
$j_0:=\mbox{max}\{j\in\nz_{\leq r}~|~k_j=j\}$, and set
$$I_1:=[j_0+1,k_r-1],
\quad\mbox{and}\quad 
I:=[j_0-(2l+1)N_{\bf a}+1,k_r+N_{\bf a}-1].$$ 
Then we have ${\bf k}\in \cM_{\snz}(I_1)\subset\cM_{\snz}(I)$. 
Hence we can consider $\Omega_I({\bf k})\in \cM_{\snz}(I)$. 
Now let us recall the definition of $M_{\Omega_I({\bf k})}({\bf a})$;
we set $I_0:=[j_0-(2l+2)N_{\bf a}+1,k_r+N_{\bf a}-1]$, and define
$$M_{\Omega_I({\bf k})}({\bf a})
=M_{\mbox{res}_{I_0}(\Omega_I({\bf k}))}({\bf a}^{I_0}).$$
Note that the explicit forms of $\mbox{res}_{I}\bigl(\Omega_I({\bf k})\bigr)$
and $\mbox{res}_{I_0}\bigl(\Omega_I({\bf k})\bigr)$ 
are given as:
\begin{align*}
\mbox{res}_{I}\bigl(\Omega_I({\bf k})\bigr)&=
\left(\widetilde{I}_1\setminus\mbox{res}_{I_1}\right)\cup[k_r+1,k_r+N_{\bf a}]
\quad\mbox{and}\\
\mbox{res}_{I_0}\bigl(\Omega_I({\bf k})\bigr)
&=[j_0-(2l+2)N_{\bf a}+1,j_0-(2l+1)N_{\bf a}]\\
&\qquad\qquad\qquad
\cup \left(\widetilde{I}_1\setminus\mbox{res}_{I_1}\right)\cup
[k_r+1,k_r+N_{\bf a}],
\end{align*}
respectively. Since $(j_0+1)-(j_0-(2l+1)N_{\bf a})>(2l+1)N_{\bf a}$, we obtain
$$M_{\Omega_I({\bf k})}({\bf a})=
M_{\mbox{res}_{I_0}(\Omega_I({\bf k}))}({\bf a}^{I_0})=
M_{\mbox{res}_{I}(\Omega_I({\bf k}))}({\bf a}^{I})\eqno{(6.3.1)}$$
by Lemma \ref{lemma:fin3}.

Let $J$ be a finite interval such that $J\supset I$. Write 
$J=[j_0-(2l+1)N_{\bf a}+1-d_1,k_r+N_{\bf a}-1+d_2]$ for some $d_1,d_2\in
\nz_{\geq 0}$, and set $J_0:=[j_0-(2l+2)N_{\bf a}+1-d_1,k_r+N_{\bf a}-1+d_2]$.
Then we have
$$M_{\Omega_J({\bf k})}({\bf a})
=M_{\mbox{res}_{J_0}(\Omega_J({\bf k}))}({\bf a}^{J_0}).$$
The explicit form of $\mbox{res}_{J_0}\bigl(\Omega_J({\bf k})\bigr)$ is
given as:
\begin{align*}
\mbox{res}_{J_0}\bigl(\Omega_J({\bf k})\bigr)
&=[j_0-(2l+2)N_{\bf a}+1-d_1,j_0-(2l+1)N_{\bf a}-d_1]\\
&\qquad\qquad\qquad
\cup \left(\widetilde{I}_1\setminus\mbox{res}_{I_1}\right)\cup
[k_r+1,k_r+N_{\bf a}+d_2].
\end{align*}
Since $(j_0+1)-(j_0-(2l+1)N_{\bf a}-d_1)>(2l+1)N_{\bf a}>N_{\bf a}$ and
$[k_r+1,k_r+N_{\bf a}]\subset [k_r+1,k_r+N_{\bf a}+d_2]$, 
we deduce that
$$M_{\Omega_J({\bf k})}({\bf a})=
M_{\mbox{res}_{J_0}(\Omega_J({\bf k}))}({\bf a}^{J_0})=
M_{\mbox{res}_{I}(\Omega_I({\bf k}))}({\bf a}^{I})$$
by part (2) of Lemma \ref{lemma:fin1}, part (1) of Lemma \ref{lemma:fin2}, 
and Lemma \ref{lemma:fin3}.
Thus, we conclude that 
$M_{\Omega_J({\bf k})}({\bf a})
=M_{\Omega_I({\bf k})}({\bf a}).$
\end{proof}
\noindent
{\it Proof of Proposition \ref{prop:e-BZ}.} Let us check 
condition (2-a) in Definition \ref{defn:BZ}. Take a finite interval $K$,
and consider a finite collection ${\bf M}({\bf a})_K=
(M_{\bf k}({\bf a}))_{{\bf k}\in\cM_{\snz}(K)}$ of integers. It suffices to prove
that ${\bf M}({\bf a})_K$ is an element of $\cBZ_{K}^{e}$ under the 
identification $\cM_{\snz}(K)\cong \cM_K^{\times}$. Since $\cM_{\snz}(K)$ is a
finite set, there exists a finite interval $I$ such that $I\supset K$ and
$M_{\bf k}({\bf a})=M_{{\bf k}^I}({\bf a}^I)$ for all ${\bf k}\in \cM_{\snz}(K)$.
Hence the condition (2-a) is satisfied by Theorem \ref{thm:finA}.

It remains to check the condition (2-b), which is clear from Lemma 
\ref{lemma:BZ-3}.\hfill$\square$
\vskip 3mm
The following corollary is easily obtained.
\begin{cor}
For each ${\bf a}\in\cB_{l-1}^{(1),ap}$, $\Phi_{\snz}({\bf a})$ is
an element of $(\cBZ_{\snz}^e)^{\sigma}$. 
\end{cor} 
\subsection{Comparison theorem}
In the previous subsection, we have constructed the map 
$\Phi_{\snz}:\cB_{l-1}^{(1),ap}\to (\cBZ_{\snz}^e)^{\sigma}$ by the assignment
${\bf a}\mapsto {\bf M}({\bf a})$. 
In this subsection, we show that
$\Phi_{\snz}$ is a morphism of crystals. 
\begin{lemma}\label{lemma:ess}
Let ${\bf a}\in \cB_{l-1}^{(1),ap}$. Then, 
\vskip 1mm
\noindent
{\rm (1)} $\mbox{\rm wt}\bigl({\bf M}({\bf a})\bigr)=\mbox{\rm wt}({\bf a})$.
\vskip 1mm
\noindent
{\rm (2)} $\heps_p^*\bigl({\bf M}({\bf a})\bigr)=\heps_p^*({\bf a})$ for
all $p\in\widehat{I}$.
\vskip 1mm
\noindent
{\rm (3)} For all $p\in\nz$, each of $\te_p^*$ and $\tf_p^*$ commutes with 
$\Phi_{\snz}$.
\end{lemma}
\begin{proof}
Recall the definition of the weight of ${\bf M}({\bf a})\in
(\cBZ_{\snz}^e)^{\sigma}$: 
$$\mbox{wt}\bigl({\bf M}({\bf a})\bigr)=\sum_{p\in \widehat{I}}
\Theta\bigl({\bf M}^*({\bf a})\bigr)_{{\bf k}(\Lambda_p)}
\widehat{\alpha}_p.$$
Here we write ${\bf M}^*({\bf a})=\bigl({\bf M}({\bf a})\bigr)^*$.
Take a sufficiently large finite interval $J$. Then, 
for all $p\in \widehat{I}$, we have
\begin{align*}
\Theta\bigl({\bf M}^*({\bf a})\bigr)_{{\bf k}(\Lambda_p)}=
M^*_{\Omega_J^c(({\bf k}(\Lambda_p))^c)}({\bf a})=
M_{(\Omega_J^c(({\bf k}(\Lambda_p))^c))^c}({\bf a})
=M_{\Omega_{J}({\bf k}(\Lambda_p))}({\bf a}),
\end{align*}
where we use (3.3.3) for the last equality. Since $J$ is sufficiently large,
we obtain the following equalities by the same argument as in the proof of
(6.3.1):
$$M_{\Omega_{J}({\bf k}(\Lambda_p))}({\bf a})=
M_{\mbox{res}_J(\Omega_{J}({\bf k}(\Lambda_p)))}({\bf a}^J)\quad
\mbox{for all }p\in \widehat{I}.$$
Write $J=[n+1,n+m]$. Note that $n+1<p<n+m$ by the construction. 
Hence we have
$$\mbox{res}_J(\Omega_{J}({\bf k}(\Lambda_p)))=[p+1,n+m+1]\in 
\cM_J^{\times}.$$
Since $J$ is sufficiently large, we may assume that $n+m-p>N_{\bf a}$. 
Therefore, we see that
$$\Theta\bigl({\bf M}^*({\bf a})\bigr)_{{\bf k}(\Lambda_p)}=
M_{[p+1,n+m+1]}({\bf a}^J)=\sum_{s\leq p}\sum_{t\geq p+1}a_{s,t}.$$
The right-hand side is exactly $r_p({\bf a})$. This proves part (1).\\

Let us prove part (2).
Since $\heps_p^*\bigl({\bf M}({\bf a})\bigr)=
\eps_p^*\bigl({\bf M}({\bf a})\bigr)$ and 
$\heps_p^*({\bf a})=\eps_p^*({\bf a})$, it suffices to show
that 
$$\eps_p^*\bigl({\bf M}({\bf a})\bigr)=\eps_p^*({\bf a})
\quad\mbox{for all }p\in\nz.
\eqno{(6.4.1)}$$
By the definitions of $\eps_p^*\bigl({\bf M}({\bf a})\bigr)$ and 
$\eps_p^*({\bf a})$, this is equivalent to:
\begin{align*}
\Theta\bigl({\bf M}^*({\bf a})\bigr)_{{\bf k}(\Lambda_p)}+
\Theta\bigl({\bf M}^*({\bf a})\bigr)_{{\bf k}(\sigma_p\Lambda_p)}
-\Theta\bigl({\bf M}^*({\bf a})\bigr)_{{\bf k}(\Lambda_{p+1})}-
\Theta\bigl({\bf M}^*({\bf a})\bigr)_{{\bf k}(\Lambda_{p-1})}
=-\eps_p({\bf a}^J)
\end{align*}
for a sufficiently large finite interval $J$.
However, this equality can be obtained by a computation similar to the above
for part (1). This proves part (2).\\

For part (3), we only give a proof of $\te_p^*\circ \Phi_{\snz}=
\Phi_{\snz}\circ \te_p^*$, since the other one follows similarly.
By (6.4.1), $\te_p^*{\bf M}({\bf a})=0$ if and only if $\te_p^*({\bf a})=0$.
Therefore, it suffices to show the commutativity under the assumption
that $\te_p^*{\bf M}({\bf a})\in (\cBZ_{\snz}^e)^{\sigma}$ (or equivalently, 
$\te_p^*{\bf a}\in \cB_{l-1}^{(1),ap}$). 

Fix a Maya diagram ${\bf k}\in\cM_{\snz}$. Then, there exists a finite interval
$J$ such that ${\bf k}\in\cM_{\snz}(J)$ and such that
the ${\bf k}$-component of $\te_p^*{\bf M}({\bf a})$ is
given by
$$\left(\te_p^*{\bf M}({\bf a})\right)_{\bf k}=
\left(\te_p^*\bigl({\bf M}({\bf a})_J\bigr)\right)_{\mbox{res}_{J}({\bf k})}.
\eqno{(6.4.2)}$$
Here, ${\bf M}({\bf a})_J
=\bigl(M_{\bf m}({\bf a})\bigr)_{{\bf m}\in \cM_{\snz}(J)}$.
Since $\cM_{\snz}(J)$ is a finite set, there exists a sufficiently
large interval $K$ such that $K\supset J$ and such that, for all 
${\bf m}\in \cM_{\snz}(J)$, 
$$M_{\bf m}({\bf a})=M_{\mbox{res}_K({\bf m})}\bigl({\bf a}^K\bigr),
\eqno{(6.4.3)}$$
$$M_{\bf m}(\te_p^*{\bf a})=M_{\mbox{res}_K({\bf m})}
\left((\te_p^*{\bf a})^K\right).\eqno{(6.4.4)}$$
Since $K$ is sufficiently large, we may assume that
$$(\te_p^*{\bf a})^K=\te_p^*\left({\bf a}^K\right).\eqno{(6.4.5)}$$

Consider a collection ${\bf M}':=\left(M_{\mbox{res}_K({\bf n})}
\bigl({\bf a}^K\bigr)\right)_{{\bf n}\in\cM_{\snz}(K)}$. It is 
the image of the Lusztig datum ${\bf a}^K$ associated to the finite interval
$K$ under the map $\Phi_{K}:\cB_K\to \cBZ_K^e$ constructed in
Subsection 2.4.1. By Theorem \ref{thm:finA} and (6.4.5), we see that
$$
\bigl(\te_p^*{\bf M}'\bigl)_{\mbox{res}_K({\bf n})}=
M_{\mbox{res}_K({\bf n})}\bigl(\te_p^*\bigl({\bf a}^K\bigr)\bigr)
=M_{\mbox{res}_K({\bf n})}\bigl(\bigl(\te_p^*{\bf a}\bigr)^K\bigr)
\eqno{(6.4.6)}
$$
for all ${\bf n}\in \cM_{\snz}(K)$. Note that the equalities above hold 
for all ${\bf m}\in \cM_{\snz}(J)$ since $\cM_{\snz}(J)\subset \cM_{\snz}(K)$. 

Set
${\bf M}'':=({\bf M}')_J=\left(M_{\mbox{res}_K({\bf m})}
\bigl({\bf a}^K\bigr)\right)_{{\bf m}\in\cM_{\snz}(J)}$. 
Then it follows that
$$\bigl(\te_p^*{\bf M}''\bigl)_{\mbox{res}_K({\bf m})}=
\bigl(\te_p^*{\bf M}'\bigl)_{\mbox{res}_K({\bf m})}\quad
\mbox{for all }{\bf m}\in\cM_{\snz}(J).$$
Therefore, by (6.4.4), (6.4.6), and these equalities, we obtain
\begin{align*}
M_{\bf m}(\te_p^*{\bf a})=\bigl(\te_p^*{\bf M}''\bigl)_{\mbox{res}_K({\bf m})}.
\end{align*} 
In particular, we have
\begin{align*}
M_{\bf k}(\te_p^*{\bf a})=\bigl(\te_p^*{\bf M}''\bigl)_{\mbox{res}_K({\bf k})}.
\end{align*} 
Also, by (6.4.3), it follows that
$${\bf M}''=\left(M_{\mbox{res}_K({\bf m})}
\bigl({\bf a}^K\bigr)\right)_{{\bf m}\in\cM_{\snz}(J)}=
\bigl({\bf M}_{\bf m}({\bf a})\bigr)_{{\bf m}\in\cM_{\snz}(J)},$$ 
from which we deduce that ${\bf M}''={\bf M}({\bf a})_J$. 
Consequently, by (6.4.2), we obtain
$$\left(\te_p^*{\bf M}({\bf a})\right)_{\bf k}=M_{\bf k}(\te_p^*{\bf a}).$$
This proves part (3).
\end{proof}

\begin{prop}\label{prop:strict}
The map $\Phi_{\snz}:\cB_{l-1}^{(1),ap}\to (\cBZ_{\snz}^e)^{\sigma}$ gives rise to 
a strict morphism of crystals form
$\left(\cB_{l-1}^{(1),ap};\mbox{\rm wt},\heps_p^*,\hvphi_p^*,\hte_p^*,\htf_p^*
\right)$ to
$\left((\cBZ_{\snz}^e)^{\sigma};
\mbox{\rm wt},\heps_p^*,\hvphi_p^*,\hte_p^*,\htf_p^*\right)$ .
\end{prop}
\begin{proof}
By part (1), (2) of Lemma \ref{lemma:ess}, it suffices to show that
$$\hte_p^*{\bf M}({\bf a})=
{\bf M}(\hte_p^*{\bf a})\quad\mbox{and}\quad
\htf_p^*{\bf M}({\bf a})=
{\bf M}(\htf_p^*{\bf a})
$$
for all $p\in\widehat{I}$ and ${\bf a}\in\cB_{l-1}^{(1),ap}$. Here it is 
understood that ${\bf M}(0)=0$. 

In the following, we only give a proof of the first equality. Note that
we may assume that $\hte_p^*{\bf M}({\bf a})\ne 0$ (or equivalently, 
$\hte_p^*{\bf a}\ne 0$). 

Let ${\bf k}\in \cM_{\snz}$. By the definitions, the ${\bf k}$-component of
$\hte_p^*{\bf M}({\bf a})$ is given as:
$$\bigl(\hte_p^*{\bf M}({\bf a})\bigr)_{\bf k}=
\left(\te_{L({\bf k}^c,p)}^*
{\bf M}({\bf a})\right)_{\bf k},$$
where $L({\bf k}^c,p)=\{q\in p+l\nz~|~q\not\in {\bf k}\mbox{ and }q+1
\in{\bf k}\}$ and $\te_{L({\bf k}^c,p)}^*=\prod_{q\in L({\bf k}^c,p)}\te_q^*$.
Therefore, by the definition of $\hte_p^*$ on $\cB_{l-1}^{(1)}$ and 
Lemma \ref{lemma:ess} (3), it suffices to verify that 
$$\bigl(\te_q^*{\bf M}({\bf a})\bigr)_{\bf k}=
\bigl({\bf M}({\bf a})\bigr)_{\bf k}\quad\mbox{for all }q\in p+l\nz\setminus
L({\bf k}^c,p).$$
Fix $q\in p+l\nz\setminus L({\bf k}^c,p)=\{q\in p+l\nz~|~q\in {\bf k}
\mbox{ or }q+1\not\in{\bf k}\}$. Then, there exists a finite interval $J$
such that $\bigl(\te_q^*{\bf M}({\bf a})\bigr)_{\bf k}=
\bigl(\te_q^*\bigl({\bf M}({\bf a})_J\bigr)\bigr)_{\mbox{res}_J(\bf k)}.$
By Corollary \ref{cor:ast-crys}, the right-hand side must be equal to
$\bigl({\bf M}({\bf a})_J\bigr)_{\mbox{res}_J(\bf k)}=
\bigl({\bf M}({\bf a})\bigr)_{\bf k}$. This proves the proposition.
\end{proof}
Now we can state the main result of this paper.
\begin{thm}\label{thm:main}
The image of $\cB_{l-1}^{(1),ap}$ under the morphism $\Phi_{\snz}$ coincides 
with $(\cBZ_{\snz}^e)^{\sigma}({\bf O}^*)$. In other words, 
$\Phi_{\snz}:
\bigl(\cB_{l-1}^{(1),ap};\mbox{\rm wt},\heps_p^*,\hvphi_p^*,\hte_p^*,\htf_p^*
\bigr)\to
\bigl((\cBZ_{\snz}^e)^{\sigma}({\bf O}^*);
\mbox{\rm wt},\heps_p^*,\hvphi_p^*,\hte_p^*,\htf_p^*\bigr)$ 
gives an isomorphism of crystals.
\end{thm}
\begin{proof}
By Theorem \ref{thm:LTV} and Proposition \ref{prop:strict}, 
$\Phi_{\snz}(\cB_{l-1}^{(1),ap})$ is a subcrystal of $(\cBZ_{\snz}^e)^{\sigma}$, 
and its crystal graph is connected. Also, 
$\Phi_{\snz}({\bf a}_{\bf 0})={\bf O}^*$. 
Therefore, $\Phi_{\snz}(\cB_{l-1}^{(1),ap})$ coincides with 
$(\cBZ_{\snz}^e)^{\sigma}({\bf O}^*)$. Namely, $\Phi_{\snz}$ is a surjective
map which preserves weights. More specifically, for $\widehat{\lambda}\in 
\widehat{Q}_-:=\oplus_{p=0}^{l-1}\nz_{\leq 0}\widehat{\alpha}_p$, set
$\bigl(\cB_{l-1}^{(1),ap}\bigr)_{\widehat{\lambda}}:=\{{\bf a}\in \cB_{l-1}^{(1),ap}
~|~\mbox{wt}({\bf a})=\widehat{\lambda}\}$ and 
$(\cBZ_{\snz}^e)^{\sigma}({\bf O}^*)_{\widehat{\lambda}}:=
\{{\bf M}\in (\cBZ_{\snz}^e)^{\sigma}({\bf O}^*)~|~\mbox{wt}({\bf M})=
\widehat{\lambda}\}$, respectively. 
Then, the restriction $\left.\Phi_{\snz}\right|_{\widehat{\lambda}}:
\bigl(\cB_{l-1}^{(1),ap}\bigr)_{\widehat{\lambda}}\to 
(\cBZ_{\snz}^e)^{\sigma}({\bf O}^*)_{\widehat{\lambda}}$ is a surjective map for
each $\widehat{\lambda}\in \widehat{Q}_-$.
Since both $\cB_{l-1}^{(1),ap}$ and 
$(\cBZ_{\snz}^e)^{\sigma}({\bf O}^*)$ are isomorphic to $B(\infty)$, the two sets
$\bigl(\cB_{l-1}^{(1),ap}\bigr)_{\widehat{\lambda}}$ and 
$(\cBZ_{\snz}^e)^{\sigma}({\bf O}^*)_{\widehat{\lambda}}$ are of the same (finite)
cardinality. Consequently, $\Phi_{\snz}$ must be a bijection. 
This proves the theorem.
\end{proof}
The following is an easy consequence of Theorem \ref{thm:main}.
\begin{cor}
The composite map $\ast\circ\Phi_{\snz}:
\cB_{l-1}^{(1),ap}
\overset{\Phi_{\snz}}{\longrightarrow}(\cBZ_{\snz}^e)^{\sigma}({\bf O}^*)
\overset{\ast}{\longrightarrow}\cBZ_{\snz}^{\sigma}({\bf O})$
gives rise to an isomorphism of crystals
$$
\left(\cB_{l-1}^{(1),ap};\mbox{\rm wt},\heps_p^*,\hvphi_p^*,\hte_p^*,\htf_p^*
\right)
\overset{\sim}{\to}
\left(\cBZ_{\snz}^{\sigma}({\bf O});
\mbox{\rm wt},\heps_p,\hvphi_p,\hte_p,\htf_p\right).$$
\end{cor}
\appendix
\section{Another explicit description}
There is another explicit description of 
$\cBZ_{\snz}^{\sigma}({\bf O})$ in terms of another 
crystal structure on $\cB_{l-1}^{(1),ap}$; namely 
$\left(\cB_{l-1}^{(1),ap};\mbox{\rm wt},\heps_p,\hvphi_p,\hte_p,\htf_p\right)$.
In the appendix, we explain it. Because all results in this appendix are
obtained by methods similar to the ones which we  
explained in this paper, we omit the details.
\subsection{}
Let $I=[n+1,n+m]$ be a finite interval in $\nz$, and 
${\bf k}=\{k_{n+1}<\cdots<k_{n+u}\}\in\cM_I^{\times}$. 
For a given Lusztig datum ${\bf a}^I=(a_{i,j}^I)_{(i,j)\in \Delta_I^+}\in\cB_I$ 
associated to
the finite interval $I$, we introduce a new collection  
${\bf M}'({\bf a}^I)=(M_{\bf k}'({\bf a}^I))_{{\bf k}\in\cM_I^{\times}}$ 
of integers by
$$M_{\bf k}'({\bf a}^I):=-\sum_{i=n+1}^{n+u}\sum_{j=k_i+1}^{n+m+1}a_{k_i,j}^I+
\mbox{min}\left\{\left.\sum_{n+1\leq p<q\leq n+u}a_{c_{p,q},c_{p,q}+(q-p)}^I
~\right|~\begin{array}{c}C=(c_{p,q})\mbox{ is }\\ 
\mbox{a ${\bf k}$-tableau.}\end{array}\right\}.$$
Define a map $\Phi_I'$ by ${\bf a}^I\mapsto {\bf M}'({\bf a}^I)$ ({\it cf}.
(2.4.1)).

By using the map $\Phi_I'$, we can construct an isomorphism from
$\cB_I$ to $\cBZ_I$ directly ({\it cf}. Corollary \ref{cor:ast-isom}).
\begin{claim}
For each ${\bf a}^I\in\cB_I$, $\Phi_I'({\bf a}^I)={\bf M}'({\bf a}^I)$ is an 
element of $\cBZ_I$. Moreover, the map $\Phi_I':\cB_I\to \cBZ_I$ is a
bijection which gives rise to an isomorphism of crystals
$$\left(\cB_I;\mbox{\rm wt},\eps_i,\vphi_i,\te_i,\tf_i\right)
\overset{\Phi_I'}{\longrightarrow}
\left(\cBZ_I;\mbox{\rm wt},\eps_i,\vphi_i,\te_i,\tf_i\right).$$
\end{claim}
\subsection{}
We discuss the infinite case. Let 
${\bf a}=(a_{i,j})_{(i,j)\in\Delta_{\snz}^+}\in\cB_{\snz}$, and set
${\bf a}^I=(a_{i,j})_{(i,j)\in \Delta_I^+}$. 
It is obvious that ${\bf a}^I$ is an element of $\cB_I$.
We can prove the following statements:
\begin{claim} 
{\rm (1)} For ${\bf a}\in\cB_{\snz}$ and ${\bf k}\in\cM^c_{\snz}$, there exists
a finite interval $I_0$ such that for every finite interval $J\supset I_0$, 
$M'_{\mbox{\rm res}^c_J({\bf k})}({\bf a}^J)=
M'_{\mbox{\rm res}^c_{I_0}({\bf k})}({\bf a}^{I_0})$. 
\vskip 1mm
\noindent
{\rm (2)} Define 
$M'_{\bf k}({\bf a}):=M'_{\mbox{\rm res}^c_{I_0}}({\bf k)}({\bf a}^{I_0})$.
Then, a collection 
${\bf M}'({\bf a}):=(M_{\bf k}'({\bf a}))_{{\bf k}\in\cM_{\snz}^c}$ is 
an element of $\cBZ_{\snz}$. In other words, we have a map
$\Phi_{\snz}':\cB_{\snz}\to\cBZ_{\snz}$ defined by ${\bf a}\mapsto 
{\bf M}'(\bf a)$.
\vskip 1mm
\noindent
{\rm (3)} The restriction of $\Phi_{\snz}'$ to $\cB_{l-1}^{(1),ap}$ gives
a bijection from $\cB_{l-1}^{(1),ap}$ to $\cBZ_{\snz}^{\sigma}({\bf O})$.
Moreover, it gives rise to an isomorphism of crystals
$$\left(\cB_{l-1}^{(1),ap};\mbox{\rm wt},\heps_p,\hvphi_p,\hte_p,\htf_p\right)
\overset{\Phi_{\snz}'}{\longrightarrow}
\left(\cBZ_{\snz}^{\sigma}({\bf O});
\mbox{\rm wt},\heps_p,\hvphi_p,\hte_p,\htf_p\right).$$
\end{claim}


\begin{thebibliography}{[MFT50]}
\bibitem[A]{A}
J. E. Anderson, 
{\it A polytope calculus for semisimple groups}, 
Duke. Math. J. {\bf 116} (2003), 567-588. 
\bibitem[BFZ]{BFZ}
A. Berenstein, S. Fomin, and A. Zelevinsky, 
{\it Parametrizations of canonical bases and totally positive matrices}, 
Adv. Math. {\bf 122} (1996), 49-149. 
\bibitem[EK]{EK}
N. Enomoto and M. Kashiwara, 
{\it Symmetric crystals for $\gtgl_{\infty}$}, 
Publ. Res. Isnt. Math. Sci. {\bf 44} (2008), 837-891. 
\bibitem[HK]{HK}
J. Hong and S. J. Kang,
{\it Introduction to Quantum Groups and Crystal Bases},
Graduate Studies in Mathematics {\bf 42} (2002), American Mathematical 
Society.
\bibitem[Kam1]{Kam1}
J. Kamnitzer,
{\it Mirkovi\'{c}-Vilonen cycles and polytopes},
Ann. of Math. {\bf 171} (2010), 245-294.
\bibitem[Kam2]{Kam2}
J. Kamnitzer,
{\it The crystal structure on the set of Mirkovi\'{c}-Vilonen polytopes},
Adv. Math. {\bf 215} (2007), 66-93. 
\bibitem[K1]{K1}
M. Kashiwara,
{\it Crystallizing the $q$-analogue of universal enveloping algebras},
Duke Math. J. {\bf 63} (1991), 465-516.
\bibitem[K2]{K2}
M. Kashiwara,
{\it Global crystal base of quantum groups},
Duke Math. J. {\bf 69} (1993), 455-485.
\bibitem[K3]{K3}
M. Kashiwara,
{\it Crystal base and Littelmann's refined Demazure character formula},
Duke Math. J. {\bf 71} (1993), 839-858.
\bibitem[K4]{K4}
M. Kashiwara,
{\it Bases Cristallines des Groupes Quantiques},
Cours Sp\'ecialis\'es Vol. {\bf 9}, Soci\'et\'e Math\'ematique de France, 
Paris, 2002.
\bibitem[KS]{KS}
M. Kashiwara and Y. Saito,
{\it Geometric construction of crystal bases},
Duke Math. J. {\bf 89} (1997), 9-36.
\bibitem[LTV]{LTV}
B. Leclerc, J-Y. Thibon, and E. Vasserot, 
{\it Zelevinsky's involution at roots of unity},
J. Reine Angew. Math. {\bf 513} (1999), 33--51.
\bibitem[L1]{L1}
G. Lusztig,
{\it Canonical bases arising from quantized universal enveloping algebras},
J. Amer. Math. Soc. {\bf 3} (1990), 447-498.
\bibitem[L2]{L3}
G. Lusztig,
{\it Quivers, perverse sheaves, and quantized universal enveloping algebras},
J. Amer. Math. Soc. {\bf 4} (1991), 365-421.
\bibitem[L3]{L2}
G. Lusztig,
{\it Affine quivers and canonical bases},
Publ. Math. IHES {\bf 76} (1992), 111-163.
\bibitem[L4]{L4}
G. Lusztig,
{\it Introduction to Quantum Groups},
Progr. Math. {\bf 110} (1993), Birkh\"auser.
\bibitem[NSS]{NSS}
S. Naito, D. Sagaki, and Y. Saito, 
{\it Toward Berenstein-Zelevinsky data in affine type $A$, I: 
Construction of affine analogs}, 
arXiv:1009.4526.
\bibitem[N]{N}
M. Noumi, 
{\it Painlev\'e Equations through Symmetry},
Translations of Mathematical Monographs {\bf 223} (2004), American 
Mathematical Society.
\bibitem[R]{R}
M. Reineke, 
{\it On the coloured graph structure of Lusztig's canonical basis}, 
Math. Ann. {\bf 307} (1997), 705-723. 
\bibitem[S]{S}
Y. Saito, 
{\it Mirkovi\'c-Vilonen polytopes and a quiver construction of
crystal basis in type $A$}, 
arXiv:1010.0086. 
\bibitem[Sav]{Sav}
A. Savage, 
{\it Geometric and combinatorial realization of crystal graphs}, 
Algebr. Represent. Theory {\bf 9} (2006), 161-199. 
\end{thebibliography}
\end{document}